\documentclass[usenames,dvipsnames]{article}

\usepackage{latexsym}
\usepackage{amsmath}
\usepackage{amsfonts}
\usepackage{amssymb}
\usepackage{psfrag} 
\usepackage{graphicx}

\usepackage{natbib}       
\usepackage{mathptmx}       
\usepackage{helvet}         
\usepackage{courier}        
\usepackage{type1cm}  
\usepackage{wrapfig}
\usepackage{float}
\usepackage{floatflt}
\usepackage{makeidx}      
\usepackage{multicol}        
\usepackage[bottom]{footmisc}
\usepackage{url}
\usepackage{tikz}
\usepackage{fixltx2e}
\usepackage{enumitem}
\usepackage{multirow}

\definecolor{AliceBlue}{rgb}{0.94,0.97,1.00}
\definecolor{AntiqueWhite1}{rgb}{1.00,0.94,0.86}
\definecolor{AntiqueWhite2}{rgb}{0.93,0.87,0.80}
\definecolor{AntiqueWhite3}{rgb}{0.80,0.75,0.69}
\definecolor{AntiqueWhite4}{rgb}{0.55,0.51,0.47}
\definecolor{AntiqueWhite}{rgb}{0.98,0.92,0.84}
\definecolor{BlanchedAlmond}{rgb}{1.00,0.92,0.80}
\definecolor{BlueViolet}{rgb}{0.54,0.17,0.89}
\definecolor{CadetBlue1}{rgb}{0.60,0.96,1.00}
\definecolor{CadetBlue2}{rgb}{0.56,0.90,0.93}
\definecolor{CadetBlue3}{rgb}{0.48,0.77,0.80}
\definecolor{CadetBlue4}{rgb}{0.33,0.53,0.55}
\definecolor{CadetBlue}{rgb}{0.37,0.62,0.63}
\definecolor{CornflowerBlue}{rgb}{0.39,0.58,0.93}
\definecolor{DarkBlue}{rgb}{0.00,0.00,0.55}
\definecolor{DarkCyan}{rgb}{0.00,0.55,0.55}
\definecolor{DarkGoldenrod1}{rgb}{1.00,0.73,0.06}
\definecolor{DarkGoldenrod2}{rgb}{0.93,0.68,0.05}
\definecolor{DarkGoldenrod3}{rgb}{0.80,0.58,0.05}
\definecolor{DarkGoldenrod4}{rgb}{0.55,0.40,0.03}
\definecolor{DarkGoldenrod}{rgb}{0.72,0.53,0.04}
\definecolor{DarkGray}{rgb}{0.66,0.66,0.66}
\definecolor{DarkGreen}{rgb}{0.00,0.39,0.00}
\definecolor{DarkGrey}{rgb}{0.66,0.66,0.66}
\definecolor{DarkKhaki}{rgb}{0.74,0.72,0.42}
\definecolor{DarkMagenta}{rgb}{0.55,0.00,0.55}
\definecolor{DarkOliveGreen1}{rgb}{0.79,1.00,0.44}
\definecolor{DarkOliveGreen2}{rgb}{0.74,0.93,0.41}
\definecolor{DarkOliveGreen3}{rgb}{0.64,0.80,0.35}
\definecolor{DarkOliveGreen4}{rgb}{0.43,0.55,0.24}
\definecolor{DarkOliveGreen}{rgb}{0.33,0.42,0.18}
\definecolor{DarkOrange1}{rgb}{1.00,0.50,0.00}
\definecolor{DarkOrange2}{rgb}{0.93,0.46,0.00}
\definecolor{DarkOrange3}{rgb}{0.80,0.40,0.00}
\definecolor{DarkOrange4}{rgb}{0.55,0.27,0.00}
\definecolor{DarkOrange}{rgb}{1.00,0.55,0.00}
\definecolor{DarkOrchid1}{rgb}{0.75,0.24,1.00}
\definecolor{DarkOrchid2}{rgb}{0.70,0.23,0.93}
\definecolor{DarkOrchid3}{rgb}{0.60,0.20,0.80}
\definecolor{DarkOrchid4}{rgb}{0.41,0.13,0.55}
\definecolor{DarkOrchid}{rgb}{0.60,0.20,0.80}
\definecolor{DarkRed}{rgb}{0.55,0.00,0.00}
\definecolor{DarkSalmon}{rgb}{0.91,0.59,0.48}
\definecolor{DarkSeaGreen1}{rgb}{0.76,1.00,0.76}
\definecolor{DarkSeaGreen2}{rgb}{0.71,0.93,0.71}
\definecolor{DarkSeaGreen3}{rgb}{0.61,0.80,0.61}
\definecolor{DarkSeaGreen4}{rgb}{0.41,0.55,0.41}
\definecolor{DarkSeaGreen}{rgb}{0.56,0.74,0.56}
\definecolor{DarkSlateBlue}{rgb}{0.28,0.24,0.55}
\definecolor{DarkSlateGray1}{rgb}{0.59,1.00,1.00}
\definecolor{DarkSlateGray2}{rgb}{0.55,0.93,0.93}
\definecolor{DarkSlateGray3}{rgb}{0.47,0.80,0.80}
\definecolor{DarkSlateGray4}{rgb}{0.32,0.55,0.55}
\definecolor{DarkSlateGray}{rgb}{0.18,0.31,0.31}
\definecolor{DarkSlateGrey}{rgb}{0.18,0.31,0.31}
\definecolor{DarkTurquoise}{rgb}{0.00,0.81,0.82}
\definecolor{DarkViolet}{rgb}{0.58,0.00,0.83}
\definecolor{DeepPink1}{rgb}{1.00,0.08,0.58}
\definecolor{DeepPink2}{rgb}{0.93,0.07,0.54}
\definecolor{DeepPink3}{rgb}{0.80,0.06,0.46}
\definecolor{DeepPink4}{rgb}{0.55,0.04,0.31}
\definecolor{DeepPink}{rgb}{1.00,0.08,0.58}
\definecolor{DeepSkyBlue1}{rgb}{0.00,0.75,1.00}
\definecolor{DeepSkyBlue2}{rgb}{0.00,0.70,0.93}
\definecolor{DeepSkyBlue3}{rgb}{0.00,0.60,0.80}
\definecolor{DeepSkyBlue4}{rgb}{0.00,0.41,0.55}
\definecolor{DeepSkyBlue}{rgb}{0.00,0.75,1.00}
\definecolor{DimGray}{rgb}{0.41,0.41,0.41}
\definecolor{DimGrey}{rgb}{0.41,0.41,0.41}
\definecolor{DodgerBlue1}{rgb}{0.12,0.56,1.00}
\definecolor{DodgerBlue2}{rgb}{0.11,0.53,0.93}
\definecolor{DodgerBlue3}{rgb}{0.09,0.45,0.80}
\definecolor{DodgerBlue4}{rgb}{0.06,0.31,0.55}
\definecolor{DodgerBlue}{rgb}{0.12,0.56,1.00}
\definecolor{FloralWhite}{rgb}{1.00,0.98,0.94}
\definecolor{ForestGreen}{rgb}{0.13,0.55,0.13}
\definecolor{GhostWhite}{rgb}{0.97,0.97,1.00}
\definecolor{GreenYellow}{rgb}{0.68,1.00,0.18}
\definecolor{HotPink1}{rgb}{1.00,0.43,0.71}
\definecolor{HotPink2}{rgb}{0.93,0.42,0.65}
\definecolor{HotPink3}{rgb}{0.80,0.38,0.56}
\definecolor{HotPink4}{rgb}{0.55,0.23,0.38}
\definecolor{HotPink}{rgb}{1.00,0.41,0.71}
\definecolor{IndianRed1}{rgb}{1.00,0.42,0.42}
\definecolor{IndianRed2}{rgb}{0.93,0.39,0.39}
\definecolor{IndianRed3}{rgb}{0.80,0.33,0.33}
\definecolor{IndianRed4}{rgb}{0.55,0.23,0.23}
\definecolor{IndianRed}{rgb}{0.80,0.36,0.36}
\definecolor{LavenderBlush1}{rgb}{1.00,0.94,0.96}
\definecolor{LavenderBlush2}{rgb}{0.93,0.88,0.90}
\definecolor{LavenderBlush3}{rgb}{0.80,0.76,0.77}
\definecolor{LavenderBlush4}{rgb}{0.55,0.51,0.53}
\definecolor{LavenderBlush}{rgb}{1.00,0.94,0.96}
\definecolor{LawnGreen}{rgb}{0.49,0.99,0.00}
\definecolor{LemonChiffon1}{rgb}{1.00,0.98,0.80}
\definecolor{LemonChiffon2}{rgb}{0.93,0.91,0.75}
\definecolor{LemonChiffon3}{rgb}{0.80,0.79,0.65}
\definecolor{LemonChiffon4}{rgb}{0.55,0.54,0.44}
\definecolor{LemonChiffon}{rgb}{1.00,0.98,0.80}
\definecolor{LightBlue1}{rgb}{0.75,0.94,1.00}
\definecolor{LightBlue2}{rgb}{0.70,0.87,0.93}
\definecolor{LightBlue3}{rgb}{0.60,0.75,0.80}
\definecolor{LightBlue4}{rgb}{0.41,0.51,0.55}
\definecolor{LightBlue}{rgb}{0.68,0.85,0.90}
\definecolor{LightCoral}{rgb}{0.94,0.50,0.50}
\definecolor{LightCyan1}{rgb}{0.88,1.00,1.00}
\definecolor{LightCyan2}{rgb}{0.82,0.93,0.93}
\definecolor{LightCyan3}{rgb}{0.71,0.80,0.80}
\definecolor{LightCyan4}{rgb}{0.48,0.55,0.55}
\definecolor{LightCyan}{rgb}{0.88,1.00,1.00}
\definecolor{LightGoldenrod1}{rgb}{1.00,0.93,0.55}
\definecolor{LightGoldenrod2}{rgb}{0.93,0.86,0.51}
\definecolor{LightGoldenrod3}{rgb}{0.80,0.75,0.44}
\definecolor{LightGoldenrod4}{rgb}{0.55,0.51,0.30}
\definecolor{LightGoldenrodYellow}{rgb}{0.98,0.98,0.82}
\definecolor{LightGoldenrod}{rgb}{0.93,0.87,0.51}
\definecolor{LightGray}{rgb}{0.83,0.83,0.83}
\definecolor{LightGreen}{rgb}{0.56,0.93,0.56}
\definecolor{LightGrey}{rgb}{0.83,0.83,0.83}
\definecolor{LightPink1}{rgb}{1.00,0.68,0.73}
\definecolor{LightPink2}{rgb}{0.93,0.64,0.68}
\definecolor{LightPink3}{rgb}{0.80,0.55,0.58}
\definecolor{LightPink4}{rgb}{0.55,0.37,0.40}
\definecolor{LightPink}{rgb}{1.00,0.71,0.76}
\definecolor{LightSalmon1}{rgb}{1.00,0.63,0.48}
\definecolor{LightSalmon2}{rgb}{0.93,0.58,0.45}
\definecolor{LightSalmon3}{rgb}{0.80,0.51,0.38}
\definecolor{LightSalmon4}{rgb}{0.55,0.34,0.26}
\definecolor{LightSalmon}{rgb}{1.00,0.63,0.48}
\definecolor{LightSeaGreen}{rgb}{0.13,0.70,0.67}
\definecolor{LightSkyBlue1}{rgb}{0.69,0.89,1.00}
\definecolor{LightSkyBlue2}{rgb}{0.64,0.83,0.93}
\definecolor{LightSkyBlue3}{rgb}{0.55,0.71,0.80}
\definecolor{LightSkyBlue4}{rgb}{0.38,0.48,0.55}
\definecolor{LightSkyBlue}{rgb}{0.53,0.81,0.98}
\definecolor{LightSlateBlue}{rgb}{0.52,0.44,1.00}
\definecolor{LightSlateGray}{rgb}{0.47,0.53,0.60}
\definecolor{LightSlateGrey}{rgb}{0.47,0.53,0.60}
\definecolor{LightSteelBlue1}{rgb}{0.79,0.88,1.00}
\definecolor{LightSteelBlue2}{rgb}{0.74,0.82,0.93}
\definecolor{LightSteelBlue3}{rgb}{0.64,0.71,0.80}
\definecolor{LightSteelBlue4}{rgb}{0.43,0.48,0.55}
\definecolor{LightSteelBlue}{rgb}{0.69,0.77,0.87}
\definecolor{LightYellow1}{rgb}{1.00,1.00,0.88}
\definecolor{LightYellow2}{rgb}{0.93,0.93,0.82}
\definecolor{LightYellow3}{rgb}{0.80,0.80,0.71}
\definecolor{LightYellow4}{rgb}{0.55,0.55,0.48}
\definecolor{LightYellow}{rgb}{1.00,1.00,0.88}
\definecolor{LimeGreen}{rgb}{0.20,0.80,0.20}
\definecolor{MediumAquamarine}{rgb}{0.40,0.80,0.67}
\definecolor{MediumBlue}{rgb}{0.00,0.00,0.80}
\definecolor{MediumOrchid1}{rgb}{0.88,0.40,1.00}
\definecolor{MediumOrchid2}{rgb}{0.82,0.37,0.93}
\definecolor{MediumOrchid3}{rgb}{0.71,0.32,0.80}
\definecolor{MediumOrchid4}{rgb}{0.48,0.22,0.55}
\definecolor{MediumOrchid}{rgb}{0.73,0.33,0.83}
\definecolor{MediumPurple1}{rgb}{0.67,0.51,1.00}
\definecolor{MediumPurple2}{rgb}{0.62,0.47,0.93}
\definecolor{MediumPurple3}{rgb}{0.54,0.41,0.80}
\definecolor{MediumPurple4}{rgb}{0.36,0.28,0.55}
\definecolor{MediumPurple}{rgb}{0.58,0.44,0.86}
\definecolor{MediumSeaGreen}{rgb}{0.24,0.70,0.44}
\definecolor{MediumSlateBlue}{rgb}{0.48,0.41,0.93}
\definecolor{MediumSpringGreen}{rgb}{0.00,0.98,0.60}
\definecolor{MediumTurquoise}{rgb}{0.28,0.82,0.80}
\definecolor{MediumVioletRed}{rgb}{0.78,0.08,0.52}
\definecolor{MidnightBlue}{rgb}{0.10,0.10,0.44}
\definecolor{MintCream}{rgb}{0.96,1.00,0.98}
\definecolor{MistyRose1}{rgb}{1.00,0.89,0.88}
\definecolor{MistyRose2}{rgb}{0.93,0.84,0.82}
\definecolor{MistyRose3}{rgb}{0.80,0.72,0.71}
\definecolor{MistyRose4}{rgb}{0.55,0.49,0.48}
\definecolor{MistyRose}{rgb}{1.00,0.89,0.88}
\definecolor{NavajoWhite1}{rgb}{1.00,0.87,0.68}
\definecolor{NavajoWhite2}{rgb}{0.93,0.81,0.63}
\definecolor{NavajoWhite3}{rgb}{0.80,0.70,0.55}
\definecolor{NavajoWhite4}{rgb}{0.55,0.47,0.37}
\definecolor{NavajoWhite}{rgb}{1.00,0.87,0.68}
\definecolor{NavyBlue}{rgb}{0.00,0.00,0.50}
\definecolor{OldLace}{rgb}{0.99,0.96,0.90}
\definecolor{OliveDrab1}{rgb}{0.75,1.00,0.24}
\definecolor{OliveDrab2}{rgb}{0.70,0.93,0.23}
\definecolor{OliveDrab3}{rgb}{0.60,0.80,0.20}
\definecolor{OliveDrab4}{rgb}{0.41,0.55,0.13}
\definecolor{OliveDrab}{rgb}{0.42,0.56,0.14}
\definecolor{OrangeRed1}{rgb}{1.00,0.27,0.00}
\definecolor{OrangeRed2}{rgb}{0.93,0.25,0.00}
\definecolor{OrangeRed3}{rgb}{0.80,0.22,0.00}
\definecolor{OrangeRed4}{rgb}{0.55,0.15,0.00}
\definecolor{OrangeRed}{rgb}{1.00,0.27,0.00}
\definecolor{PaleGoldenrod}{rgb}{0.93,0.91,0.67}
\definecolor{PaleGreen1}{rgb}{0.60,1.00,0.60}
\definecolor{PaleGreen2}{rgb}{0.56,0.93,0.56}
\definecolor{PaleGreen3}{rgb}{0.49,0.80,0.49}
\definecolor{PaleGreen4}{rgb}{0.33,0.55,0.33}
\definecolor{PaleGreen}{rgb}{0.60,0.98,0.60}
\definecolor{PaleTurquoise1}{rgb}{0.73,1.00,1.00}
\definecolor{PaleTurquoise2}{rgb}{0.68,0.93,0.93}
\definecolor{PaleTurquoise3}{rgb}{0.59,0.80,0.80}
\definecolor{PaleTurquoise4}{rgb}{0.40,0.55,0.55}
\definecolor{PaleTurquoise}{rgb}{0.69,0.93,0.93}
\definecolor{PaleVioletRed1}{rgb}{1.00,0.51,0.67}
\definecolor{PaleVioletRed2}{rgb}{0.93,0.47,0.62}
\definecolor{PaleVioletRed3}{rgb}{0.80,0.41,0.54}
\definecolor{PaleVioletRed4}{rgb}{0.55,0.28,0.36}
\definecolor{PaleVioletRed}{rgb}{0.86,0.44,0.58}
\definecolor{PapayaWhip}{rgb}{1.00,0.94,0.84}
\definecolor{PeachPuff1}{rgb}{1.00,0.85,0.73}
\definecolor{PeachPuff2}{rgb}{0.93,0.80,0.68}
\definecolor{PeachPuff3}{rgb}{0.80,0.69,0.58}
\definecolor{PeachPuff4}{rgb}{0.55,0.47,0.40}
\definecolor{PeachPuff}{rgb}{1.00,0.85,0.73}
\definecolor{PowderBlue}{rgb}{0.69,0.88,0.90}
\definecolor{RosyBrown1}{rgb}{1.00,0.76,0.76}
\definecolor{RosyBrown2}{rgb}{0.93,0.71,0.71}
\definecolor{RosyBrown3}{rgb}{0.80,0.61,0.61}
\definecolor{RosyBrown4}{rgb}{0.55,0.41,0.41}
\definecolor{RosyBrown}{rgb}{0.74,0.56,0.56}
\definecolor{RoyalBlue1}{rgb}{0.28,0.46,1.00}
\definecolor{RoyalBlue2}{rgb}{0.26,0.43,0.93}
\definecolor{RoyalBlue3}{rgb}{0.23,0.37,0.80}
\definecolor{RoyalBlue4}{rgb}{0.15,0.25,0.55}
\definecolor{RoyalBlue}{rgb}{0.25,0.41,0.88}
\definecolor{SaddleBrown}{rgb}{0.55,0.27,0.07}
\definecolor{SandyBrown}{rgb}{0.96,0.64,0.38}
\definecolor{SeaGreen1}{rgb}{0.33,1.00,0.62}
\definecolor{SeaGreen2}{rgb}{0.31,0.93,0.58}
\definecolor{SeaGreen3}{rgb}{0.26,0.80,0.50}
\definecolor{SeaGreen4}{rgb}{0.18,0.55,0.34}
\definecolor{SeaGreen}{rgb}{0.18,0.55,0.34}
\definecolor{SkyBlue1}{rgb}{0.53,0.81,1.00}
\definecolor{SkyBlue2}{rgb}{0.49,0.75,0.93}
\definecolor{SkyBlue3}{rgb}{0.42,0.65,0.80}
\definecolor{SkyBlue4}{rgb}{0.29,0.44,0.55}
\definecolor{SkyBlue}{rgb}{0.53,0.81,0.92}
\definecolor{SlateBlue1}{rgb}{0.51,0.44,1.00}
\definecolor{SlateBlue2}{rgb}{0.48,0.40,0.93}
\definecolor{SlateBlue3}{rgb}{0.41,0.35,0.80}
\definecolor{SlateBlue4}{rgb}{0.28,0.24,0.55}
\definecolor{SlateBlue}{rgb}{0.42,0.35,0.80}
\definecolor{SlateGray1}{rgb}{0.78,0.89,1.00}
\definecolor{SlateGray2}{rgb}{0.73,0.83,0.93}
\definecolor{SlateGray3}{rgb}{0.62,0.71,0.80}
\definecolor{SlateGray4}{rgb}{0.42,0.48,0.55}
\definecolor{SlateGray}{rgb}{0.44,0.50,0.56}
\definecolor{SlateGrey}{rgb}{0.44,0.50,0.56}
\definecolor{SpringGreen1}{rgb}{0.00,1.00,0.50}
\definecolor{SpringGreen2}{rgb}{0.00,0.93,0.46}
\definecolor{SpringGreen3}{rgb}{0.00,0.80,0.40}
\definecolor{SpringGreen4}{rgb}{0.00,0.55,0.27}
\definecolor{SpringGreen}{rgb}{0.00,1.00,0.50}
\definecolor{SteelBlue1}{rgb}{0.39,0.72,1.00}
\definecolor{SteelBlue2}{rgb}{0.36,0.67,0.93}
\definecolor{SteelBlue3}{rgb}{0.31,0.58,0.80}
\definecolor{SteelBlue4}{rgb}{0.21,0.39,0.55}
\definecolor{SteelBlue}{rgb}{0.27,0.51,0.71}
\definecolor{VioletRed1}{rgb}{1.00,0.24,0.59}
\definecolor{VioletRed2}{rgb}{0.93,0.23,0.55}
\definecolor{VioletRed3}{rgb}{0.80,0.20,0.47}
\definecolor{VioletRed4}{rgb}{0.55,0.13,0.32}
\definecolor{VioletRed}{rgb}{0.82,0.13,0.56}
\definecolor{WhiteSmoke}{rgb}{0.96,0.96,0.96}
\definecolor{YellowGreen}{rgb}{0.60,0.80,0.20}
\definecolor{aliceblue}{rgb}{0.94,0.97,1.00}
\definecolor{antiquewhite}{rgb}{0.98,0.92,0.84}
\definecolor{aquamarine1}{rgb}{0.50,1.00,0.83}
\definecolor{aquamarine2}{rgb}{0.46,0.93,0.78}
\definecolor{aquamarine3}{rgb}{0.40,0.80,0.67}
\definecolor{aquamarine4}{rgb}{0.27,0.55,0.45}
\definecolor{aquamarine}{rgb}{0.50,1.00,0.83}
\definecolor{azure1}{rgb}{0.94,1.00,1.00}
\definecolor{azure2}{rgb}{0.88,0.93,0.93}
\definecolor{azure3}{rgb}{0.76,0.80,0.80}
\definecolor{azure4}{rgb}{0.51,0.55,0.55}
\definecolor{azure}{rgb}{0.94,1.00,1.00}
\definecolor{beige}{rgb}{0.96,0.96,0.86}
\definecolor{bisque1}{rgb}{1.00,0.89,0.77}
\definecolor{bisque2}{rgb}{0.93,0.84,0.72}
\definecolor{bisque3}{rgb}{0.80,0.72,0.62}
\definecolor{bisque4}{rgb}{0.55,0.49,0.42}
\definecolor{bisque}{rgb}{1.00,0.89,0.77}
\definecolor{black}{rgb}{0.00,0.00,0.00}
\definecolor{blanchedalmond}{rgb}{1.00,0.92,0.80}
\definecolor{blue1}{rgb}{0.00,0.00,1.00}
\definecolor{blue2}{rgb}{0.00,0.00,0.93}
\definecolor{blue3}{rgb}{0.00,0.00,0.80}
\definecolor{blue4}{rgb}{0.00,0.00,0.55}
\definecolor{blueviolet}{rgb}{0.54,0.17,0.89}
\definecolor{blue}{rgb}{0.00,0.00,1.00}
\definecolor{brown1}{rgb}{1.00,0.25,0.25}
\definecolor{brown2}{rgb}{0.93,0.23,0.23}
\definecolor{brown3}{rgb}{0.80,0.20,0.20}
\definecolor{brown4}{rgb}{0.55,0.14,0.14}
\definecolor{brown}{rgb}{0.65,0.16,0.16}
\definecolor{burlywood1}{rgb}{1.00,0.83,0.61}
\definecolor{burlywood2}{rgb}{0.93,0.77,0.57}
\definecolor{burlywood3}{rgb}{0.80,0.67,0.49}
\definecolor{burlywood4}{rgb}{0.55,0.45,0.33}
\definecolor{burlywood}{rgb}{0.87,0.72,0.53}
\definecolor{cadetblue}{rgb}{0.37,0.62,0.63}
\definecolor{chartreuse1}{rgb}{0.50,1.00,0.00}
\definecolor{chartreuse2}{rgb}{0.46,0.93,0.00}
\definecolor{chartreuse3}{rgb}{0.40,0.80,0.00}
\definecolor{chartreuse4}{rgb}{0.27,0.55,0.00}
\definecolor{chartreuse}{rgb}{0.50,1.00,0.00}
\definecolor{chocolate1}{rgb}{1.00,0.50,0.14}
\definecolor{chocolate2}{rgb}{0.93,0.46,0.13}
\definecolor{chocolate3}{rgb}{0.80,0.40,0.11}
\definecolor{chocolate4}{rgb}{0.55,0.27,0.07}
\definecolor{chocolate}{rgb}{0.82,0.41,0.12}
\definecolor{coral1}{rgb}{1.00,0.45,0.34}
\definecolor{coral2}{rgb}{0.93,0.42,0.31}
\definecolor{coral3}{rgb}{0.80,0.36,0.27}
\definecolor{coral4}{rgb}{0.55,0.24,0.18}
\definecolor{coral}{rgb}{1.00,0.50,0.31}
\definecolor{cornflowerblue}{rgb}{0.39,0.58,0.93}
\definecolor{cornsilk1}{rgb}{1.00,0.97,0.86}
\definecolor{cornsilk2}{rgb}{0.93,0.91,0.80}
\definecolor{cornsilk3}{rgb}{0.80,0.78,0.69}
\definecolor{cornsilk4}{rgb}{0.55,0.53,0.47}
\definecolor{cornsilk}{rgb}{1.00,0.97,0.86}
\definecolor{cyan1}{rgb}{0.00,1.00,1.00}
\definecolor{cyan2}{rgb}{0.00,0.93,0.93}
\definecolor{cyan3}{rgb}{0.00,0.80,0.80}
\definecolor{cyan4}{rgb}{0.00,0.55,0.55}
\definecolor{cyan}{rgb}{0.00,1.00,1.00}
\definecolor{darkblue}{rgb}{0.00,0.00,0.55}
\definecolor{darkcyan}{rgb}{0.00,0.55,0.55}
\definecolor{darkgoldenrod}{rgb}{0.72,0.53,0.04}
\definecolor{darkgray}{rgb}{0.66,0.66,0.66}
\definecolor{darkgreen}{rgb}{0.00,0.39,0.00}
\definecolor{darkgrey}{rgb}{0.66,0.66,0.66}
\definecolor{darkkhaki}{rgb}{0.74,0.72,0.42}
\definecolor{darkmagenta}{rgb}{0.55,0.00,0.55}
\definecolor{darkolive}{rgb}{0.33,0.42,0.18}
\definecolor{darkorange}{rgb}{1.00,0.55,0.00}
\definecolor{darkorchid}{rgb}{0.60,0.20,0.80}
\definecolor{darkred}{rgb}{0.55,0.00,0.00}
\definecolor{darksalmon}{rgb}{0.91,0.59,0.48}
\definecolor{darksea}{rgb}{0.56,0.74,0.56}
\definecolor{darkslate}{rgb}{0.18,0.31,0.31}
\definecolor{darkslate}{rgb}{0.18,0.31,0.31}
\definecolor{darkslate}{rgb}{0.28,0.24,0.55}
\definecolor{darkturquoise}{rgb}{0.00,0.81,0.82}
\definecolor{darkviolet}{rgb}{0.58,0.00,0.83}
\definecolor{deeppink}{rgb}{1.00,0.08,0.58}
\definecolor{deepsky}{rgb}{0.00,0.75,1.00}
\definecolor{dimgray}{rgb}{0.41,0.41,0.41}
\definecolor{dimgrey}{rgb}{0.41,0.41,0.41}
\definecolor{dodgerblue}{rgb}{0.12,0.56,1.00}
\definecolor{firebrick1}{rgb}{1.00,0.19,0.19}
\definecolor{firebrick2}{rgb}{0.93,0.17,0.17}
\definecolor{firebrick3}{rgb}{0.80,0.15,0.15}
\definecolor{firebrick4}{rgb}{0.55,0.10,0.10}
\definecolor{firebrick}{rgb}{0.70,0.13,0.13}
\definecolor{floralwhite}{rgb}{1.00,0.98,0.94}
\definecolor{forestgreen}{rgb}{0.13,0.55,0.13}
\definecolor{gainsboro}{rgb}{0.86,0.86,0.86}
\definecolor{ghostwhite}{rgb}{0.97,0.97,1.00}
\definecolor{gold1}{rgb}{1.00,0.84,0.00}
\definecolor{gold2}{rgb}{0.93,0.79,0.00}
\definecolor{gold3}{rgb}{0.80,0.68,0.00}
\definecolor{gold4}{rgb}{0.55,0.46,0.00}
\definecolor{goldenrod1}{rgb}{1.00,0.76,0.15}
\definecolor{goldenrod2}{rgb}{0.93,0.71,0.13}
\definecolor{goldenrod3}{rgb}{0.80,0.61,0.11}
\definecolor{goldenrod4}{rgb}{0.55,0.41,0.08}
\definecolor{goldenrod}{rgb}{0.85,0.65,0.13}
\definecolor{gold}{rgb}{1.00,0.84,0.00}
\definecolor{gray0}{rgb}{0.00,0.00,0.00}
\definecolor{gray100}{rgb}{1.00,1.00,1.00}
\definecolor{gray10}{rgb}{0.10,0.10,0.10}
\definecolor{gray11}{rgb}{0.11,0.11,0.11}
\definecolor{gray12}{rgb}{0.12,0.12,0.12}
\definecolor{gray13}{rgb}{0.13,0.13,0.13}
\definecolor{gray14}{rgb}{0.14,0.14,0.14}
\definecolor{gray15}{rgb}{0.15,0.15,0.15}
\definecolor{gray16}{rgb}{0.16,0.16,0.16}
\definecolor{gray17}{rgb}{0.17,0.17,0.17}
\definecolor{gray18}{rgb}{0.18,0.18,0.18}
\definecolor{gray19}{rgb}{0.19,0.19,0.19}
\definecolor{gray1}{rgb}{0.01,0.01,0.01}
\definecolor{gray20}{rgb}{0.20,0.20,0.20}
\definecolor{gray21}{rgb}{0.21,0.21,0.21}
\definecolor{gray22}{rgb}{0.22,0.22,0.22}
\definecolor{gray23}{rgb}{0.23,0.23,0.23}
\definecolor{gray24}{rgb}{0.24,0.24,0.24}
\definecolor{gray25}{rgb}{0.25,0.25,0.25}
\definecolor{gray26}{rgb}{0.26,0.26,0.26}
\definecolor{gray27}{rgb}{0.27,0.27,0.27}
\definecolor{gray28}{rgb}{0.28,0.28,0.28}
\definecolor{gray29}{rgb}{0.29,0.29,0.29}
\definecolor{gray2}{rgb}{0.02,0.02,0.02}
\definecolor{gray30}{rgb}{0.30,0.30,0.30}
\definecolor{gray31}{rgb}{0.31,0.31,0.31}
\definecolor{gray32}{rgb}{0.32,0.32,0.32}
\definecolor{gray33}{rgb}{0.33,0.33,0.33}
\definecolor{gray34}{rgb}{0.34,0.34,0.34}
\definecolor{gray35}{rgb}{0.35,0.35,0.35}
\definecolor{gray36}{rgb}{0.36,0.36,0.36}
\definecolor{gray37}{rgb}{0.37,0.37,0.37}
\definecolor{gray38}{rgb}{0.38,0.38,0.38}
\definecolor{gray39}{rgb}{0.39,0.39,0.39}
\definecolor{gray3}{rgb}{0.03,0.03,0.03}
\definecolor{gray40}{rgb}{0.40,0.40,0.40}
\definecolor{gray41}{rgb}{0.41,0.41,0.41}
\definecolor{gray42}{rgb}{0.42,0.42,0.42}
\definecolor{gray43}{rgb}{0.43,0.43,0.43}
\definecolor{gray44}{rgb}{0.44,0.44,0.44}
\definecolor{gray45}{rgb}{0.45,0.45,0.45}
\definecolor{gray46}{rgb}{0.46,0.46,0.46}
\definecolor{gray47}{rgb}{0.47,0.47,0.47}
\definecolor{gray48}{rgb}{0.48,0.48,0.48}
\definecolor{gray49}{rgb}{0.49,0.49,0.49}
\definecolor{gray4}{rgb}{0.04,0.04,0.04}
\definecolor{gray50}{rgb}{0.50,0.50,0.50}
\definecolor{gray51}{rgb}{0.51,0.51,0.51}
\definecolor{gray52}{rgb}{0.52,0.52,0.52}
\definecolor{gray53}{rgb}{0.53,0.53,0.53}
\definecolor{gray54}{rgb}{0.54,0.54,0.54}
\definecolor{gray55}{rgb}{0.55,0.55,0.55}
\definecolor{gray56}{rgb}{0.56,0.56,0.56}
\definecolor{gray57}{rgb}{0.57,0.57,0.57}
\definecolor{gray58}{rgb}{0.58,0.58,0.58}
\definecolor{gray59}{rgb}{0.59,0.59,0.59}
\definecolor{gray5}{rgb}{0.05,0.05,0.05}
\definecolor{gray60}{rgb}{0.60,0.60,0.60}
\definecolor{gray61}{rgb}{0.61,0.61,0.61}
\definecolor{gray62}{rgb}{0.62,0.62,0.62}
\definecolor{gray63}{rgb}{0.63,0.63,0.63}
\definecolor{gray64}{rgb}{0.64,0.64,0.64}
\definecolor{gray65}{rgb}{0.65,0.65,0.65}
\definecolor{gray66}{rgb}{0.66,0.66,0.66}
\definecolor{gray67}{rgb}{0.67,0.67,0.67}
\definecolor{gray68}{rgb}{0.68,0.68,0.68}
\definecolor{gray69}{rgb}{0.69,0.69,0.69}
\definecolor{gray6}{rgb}{0.06,0.06,0.06}
\definecolor{gray70}{rgb}{0.70,0.70,0.70}
\definecolor{gray71}{rgb}{0.71,0.71,0.71}
\definecolor{gray72}{rgb}{0.72,0.72,0.72}
\definecolor{gray73}{rgb}{0.73,0.73,0.73}
\definecolor{gray74}{rgb}{0.74,0.74,0.74}
\definecolor{gray75}{rgb}{0.75,0.75,0.75}
\definecolor{gray76}{rgb}{0.76,0.76,0.76}
\definecolor{gray77}{rgb}{0.77,0.77,0.77}
\definecolor{gray78}{rgb}{0.78,0.78,0.78}
\definecolor{gray79}{rgb}{0.79,0.79,0.79}
\definecolor{gray7}{rgb}{0.07,0.07,0.07}
\definecolor{gray80}{rgb}{0.80,0.80,0.80}
\definecolor{gray81}{rgb}{0.81,0.81,0.81}
\definecolor{gray82}{rgb}{0.82,0.82,0.82}
\definecolor{gray83}{rgb}{0.83,0.83,0.83}
\definecolor{gray84}{rgb}{0.84,0.84,0.84}
\definecolor{gray85}{rgb}{0.85,0.85,0.85}
\definecolor{gray86}{rgb}{0.86,0.86,0.86}
\definecolor{gray87}{rgb}{0.87,0.87,0.87}
\definecolor{gray88}{rgb}{0.88,0.88,0.88}
\definecolor{gray89}{rgb}{0.89,0.89,0.89}
\definecolor{gray8}{rgb}{0.08,0.08,0.08}
\definecolor{gray90}{rgb}{0.90,0.90,0.90}
\definecolor{gray91}{rgb}{0.91,0.91,0.91}
\definecolor{gray92}{rgb}{0.92,0.92,0.92}
\definecolor{gray93}{rgb}{0.93,0.93,0.93}
\definecolor{gray94}{rgb}{0.94,0.94,0.94}
\definecolor{gray95}{rgb}{0.95,0.95,0.95}
\definecolor{gray96}{rgb}{0.96,0.96,0.96}
\definecolor{gray97}{rgb}{0.97,0.97,0.97}
\definecolor{gray98}{rgb}{0.98,0.98,0.98}
\definecolor{gray99}{rgb}{0.99,0.99,0.99}
\definecolor{gray9}{rgb}{0.09,0.09,0.09}
\definecolor{gray}{rgb}{0.75,0.75,0.75}
\definecolor{green1}{rgb}{0.00,1.00,0.00}
\definecolor{green2}{rgb}{0.00,0.93,0.00}
\definecolor{green3}{rgb}{0.00,0.80,0.00}
\definecolor{green4}{rgb}{0.00,0.55,0.00}
\definecolor{greenyellow}{rgb}{0.68,1.00,0.18}
\definecolor{green}{rgb}{0.00,1.00,0.00}
\definecolor{grey0}{rgb}{0.00,0.00,0.00}
\definecolor{grey100}{rgb}{1.00,1.00,1.00}
\definecolor{grey10}{rgb}{0.10,0.10,0.10}
\definecolor{grey11}{rgb}{0.11,0.11,0.11}
\definecolor{grey12}{rgb}{0.12,0.12,0.12}
\definecolor{grey13}{rgb}{0.13,0.13,0.13}
\definecolor{grey14}{rgb}{0.14,0.14,0.14}
\definecolor{grey15}{rgb}{0.15,0.15,0.15}
\definecolor{grey16}{rgb}{0.16,0.16,0.16}
\definecolor{grey17}{rgb}{0.17,0.17,0.17}
\definecolor{grey18}{rgb}{0.18,0.18,0.18}
\definecolor{grey19}{rgb}{0.19,0.19,0.19}
\definecolor{grey1}{rgb}{0.01,0.01,0.01}
\definecolor{grey20}{rgb}{0.20,0.20,0.20}
\definecolor{grey21}{rgb}{0.21,0.21,0.21}
\definecolor{grey22}{rgb}{0.22,0.22,0.22}
\definecolor{grey23}{rgb}{0.23,0.23,0.23}
\definecolor{grey24}{rgb}{0.24,0.24,0.24}
\definecolor{grey25}{rgb}{0.25,0.25,0.25}
\definecolor{grey26}{rgb}{0.26,0.26,0.26}
\definecolor{grey27}{rgb}{0.27,0.27,0.27}
\definecolor{grey28}{rgb}{0.28,0.28,0.28}
\definecolor{grey29}{rgb}{0.29,0.29,0.29}
\definecolor{grey2}{rgb}{0.02,0.02,0.02}
\definecolor{grey30}{rgb}{0.30,0.30,0.30}
\definecolor{grey31}{rgb}{0.31,0.31,0.31}
\definecolor{grey32}{rgb}{0.32,0.32,0.32}
\definecolor{grey33}{rgb}{0.33,0.33,0.33}
\definecolor{grey34}{rgb}{0.34,0.34,0.34}
\definecolor{grey35}{rgb}{0.35,0.35,0.35}
\definecolor{grey36}{rgb}{0.36,0.36,0.36}
\definecolor{grey37}{rgb}{0.37,0.37,0.37}
\definecolor{grey38}{rgb}{0.38,0.38,0.38}
\definecolor{grey39}{rgb}{0.39,0.39,0.39}
\definecolor{grey3}{rgb}{0.03,0.03,0.03}
\definecolor{grey40}{rgb}{0.40,0.40,0.40}
\definecolor{grey41}{rgb}{0.41,0.41,0.41}
\definecolor{grey42}{rgb}{0.42,0.42,0.42}
\definecolor{grey43}{rgb}{0.43,0.43,0.43}
\definecolor{grey44}{rgb}{0.44,0.44,0.44}
\definecolor{grey45}{rgb}{0.45,0.45,0.45}
\definecolor{grey46}{rgb}{0.46,0.46,0.46}
\definecolor{grey47}{rgb}{0.47,0.47,0.47}
\definecolor{grey48}{rgb}{0.48,0.48,0.48}
\definecolor{grey49}{rgb}{0.49,0.49,0.49}
\definecolor{grey4}{rgb}{0.04,0.04,0.04}
\definecolor{grey50}{rgb}{0.50,0.50,0.50}
\definecolor{grey51}{rgb}{0.51,0.51,0.51}
\definecolor{grey52}{rgb}{0.52,0.52,0.52}
\definecolor{grey53}{rgb}{0.53,0.53,0.53}
\definecolor{grey54}{rgb}{0.54,0.54,0.54}
\definecolor{grey55}{rgb}{0.55,0.55,0.55}
\definecolor{grey56}{rgb}{0.56,0.56,0.56}
\definecolor{grey57}{rgb}{0.57,0.57,0.57}
\definecolor{grey58}{rgb}{0.58,0.58,0.58}
\definecolor{grey59}{rgb}{0.59,0.59,0.59}
\definecolor{grey5}{rgb}{0.05,0.05,0.05}
\definecolor{grey60}{rgb}{0.60,0.60,0.60}
\definecolor{grey61}{rgb}{0.61,0.61,0.61}
\definecolor{grey62}{rgb}{0.62,0.62,0.62}
\definecolor{grey63}{rgb}{0.63,0.63,0.63}
\definecolor{grey64}{rgb}{0.64,0.64,0.64}
\definecolor{grey65}{rgb}{0.65,0.65,0.65}
\definecolor{grey66}{rgb}{0.66,0.66,0.66}
\definecolor{grey67}{rgb}{0.67,0.67,0.67}
\definecolor{grey68}{rgb}{0.68,0.68,0.68}
\definecolor{grey69}{rgb}{0.69,0.69,0.69}
\definecolor{grey6}{rgb}{0.06,0.06,0.06}
\definecolor{grey70}{rgb}{0.70,0.70,0.70}
\definecolor{grey71}{rgb}{0.71,0.71,0.71}
\definecolor{grey72}{rgb}{0.72,0.72,0.72}
\definecolor{grey73}{rgb}{0.73,0.73,0.73}
\definecolor{grey74}{rgb}{0.74,0.74,0.74}
\definecolor{grey75}{rgb}{0.75,0.75,0.75}
\definecolor{grey76}{rgb}{0.76,0.76,0.76}
\definecolor{grey77}{rgb}{0.77,0.77,0.77}
\definecolor{grey78}{rgb}{0.78,0.78,0.78}
\definecolor{grey79}{rgb}{0.79,0.79,0.79}
\definecolor{grey7}{rgb}{0.07,0.07,0.07}
\definecolor{grey80}{rgb}{0.80,0.80,0.80}
\definecolor{grey81}{rgb}{0.81,0.81,0.81}
\definecolor{grey82}{rgb}{0.82,0.82,0.82}
\definecolor{grey83}{rgb}{0.83,0.83,0.83}
\definecolor{grey84}{rgb}{0.84,0.84,0.84}
\definecolor{grey85}{rgb}{0.85,0.85,0.85}
\definecolor{grey86}{rgb}{0.86,0.86,0.86}
\definecolor{grey87}{rgb}{0.87,0.87,0.87}
\definecolor{grey88}{rgb}{0.88,0.88,0.88}
\definecolor{grey89}{rgb}{0.89,0.89,0.89}
\definecolor{grey8}{rgb}{0.08,0.08,0.08}
\definecolor{grey90}{rgb}{0.90,0.90,0.90}
\definecolor{grey91}{rgb}{0.91,0.91,0.91}
\definecolor{grey92}{rgb}{0.92,0.92,0.92}
\definecolor{grey93}{rgb}{0.93,0.93,0.93}
\definecolor{grey94}{rgb}{0.94,0.94,0.94}
\definecolor{grey95}{rgb}{0.95,0.95,0.95}
\definecolor{grey96}{rgb}{0.96,0.96,0.96}
\definecolor{grey97}{rgb}{0.97,0.97,0.97}
\definecolor{grey98}{rgb}{0.98,0.98,0.98}
\definecolor{grey99}{rgb}{0.99,0.99,0.99}
\definecolor{grey9}{rgb}{0.09,0.09,0.09}
\definecolor{grey}{rgb}{0.75,0.75,0.75}
\definecolor{honeydew1}{rgb}{0.94,1.00,0.94}
\definecolor{honeydew2}{rgb}{0.88,0.93,0.88}
\definecolor{honeydew3}{rgb}{0.76,0.80,0.76}
\definecolor{honeydew4}{rgb}{0.51,0.55,0.51}
\definecolor{honeydew}{rgb}{0.94,1.00,0.94}
\definecolor{hotpink}{rgb}{1.00,0.41,0.71}
\definecolor{indianred}{rgb}{0.80,0.36,0.36}
\definecolor{ivory1}{rgb}{1.00,1.00,0.94}
\definecolor{ivory2}{rgb}{0.93,0.93,0.88}
\definecolor{ivory3}{rgb}{0.80,0.80,0.76}
\definecolor{ivory4}{rgb}{0.55,0.55,0.51}
\definecolor{ivory}{rgb}{1.00,1.00,0.94}
\definecolor{khaki1}{rgb}{1.00,0.96,0.56}
\definecolor{khaki2}{rgb}{0.93,0.90,0.52}
\definecolor{khaki3}{rgb}{0.80,0.78,0.45}
\definecolor{khaki4}{rgb}{0.55,0.53,0.31}
\definecolor{khaki}{rgb}{0.94,0.90,0.55}
\definecolor{lavenderblush}{rgb}{1.00,0.94,0.96}
\definecolor{lavender}{rgb}{0.90,0.90,0.98}
\definecolor{lawngreen}{rgb}{0.49,0.99,0.00}
\definecolor{lemonchiffon}{rgb}{1.00,0.98,0.80}
\definecolor{lightblue}{rgb}{0.68,0.85,0.90}
\definecolor{lightcoral}{rgb}{0.94,0.50,0.50}
\definecolor{lightcyan}{rgb}{0.88,1.00,1.00}
\definecolor{lightgoldenrod}{rgb}{0.93,0.87,0.51}
\definecolor{lightgoldenrod}{rgb}{0.98,0.98,0.82}
\definecolor{lightgray}{rgb}{0.83,0.83,0.83}
\definecolor{lightgreen}{rgb}{0.56,0.93,0.56}
\definecolor{lightgrey}{rgb}{0.83,0.83,0.83}
\definecolor{lightpink}{rgb}{1.00,0.71,0.76}
\definecolor{lightsalmon}{rgb}{1.00,0.63,0.48}
\definecolor{lightsea}{rgb}{0.13,0.70,0.67}
\definecolor{lightsky}{rgb}{0.53,0.81,0.98}
\definecolor{lightslate}{rgb}{0.47,0.53,0.60}
\definecolor{lightslate}{rgb}{0.47,0.53,0.60}
\definecolor{lightslate}{rgb}{0.52,0.44,1.00}
\definecolor{lightsteel}{rgb}{0.69,0.77,0.87}
\definecolor{lightyellow}{rgb}{1.00,1.00,0.88}
\definecolor{limegreen}{rgb}{0.20,0.80,0.20}
\definecolor{linen}{rgb}{0.98,0.94,0.90}
\definecolor{magenta1}{rgb}{1.00,0.00,1.00}
\definecolor{magenta2}{rgb}{0.93,0.00,0.93}
\definecolor{magenta3}{rgb}{0.80,0.00,0.80}
\definecolor{magenta4}{rgb}{0.55,0.00,0.55}
\definecolor{magenta}{rgb}{1.00,0.00,1.00}
\definecolor{maroon1}{rgb}{1.00,0.20,0.70}
\definecolor{maroon2}{rgb}{0.93,0.19,0.65}
\definecolor{maroon3}{rgb}{0.80,0.16,0.56}
\definecolor{maroon4}{rgb}{0.55,0.11,0.38}
\definecolor{maroon}{rgb}{0.69,0.19,0.38}
\definecolor{mediumaquamarine}{rgb}{0.40,0.80,0.67}
\definecolor{mediumblue}{rgb}{0.00,0.00,0.80}
\definecolor{mediumorchid}{rgb}{0.73,0.33,0.83}
\definecolor{mediumpurple}{rgb}{0.58,0.44,0.86}
\definecolor{mediumsea}{rgb}{0.24,0.70,0.44}
\definecolor{mediumslate}{rgb}{0.48,0.41,0.93}
\definecolor{mediumspring}{rgb}{0.00,0.98,0.60}
\definecolor{mediumturquoise}{rgb}{0.28,0.82,0.80}
\definecolor{mediumviolet}{rgb}{0.78,0.08,0.52}
\definecolor{midnightblue}{rgb}{0.10,0.10,0.44}
\definecolor{mintcream}{rgb}{0.96,1.00,0.98}
\definecolor{mistyrose}{rgb}{1.00,0.89,0.88}
\definecolor{moccasin}{rgb}{1.00,0.89,0.71}
\definecolor{navajowhite}{rgb}{1.00,0.87,0.68}
\definecolor{navyblue}{rgb}{0.00,0.00,0.50}
\definecolor{navy}{rgb}{0.00,0.00,0.50}
\definecolor{oldlace}{rgb}{0.99,0.96,0.90}
\definecolor{olivedrab}{rgb}{0.42,0.56,0.14}
\definecolor{orange1}{rgb}{1.00,0.65,0.00}
\definecolor{orange2}{rgb}{0.93,0.60,0.00}
\definecolor{orange3}{rgb}{0.80,0.52,0.00}
\definecolor{orange4}{rgb}{0.55,0.35,0.00}
\definecolor{orangered}{rgb}{1.00,0.27,0.00}
\definecolor{orange}{rgb}{1.00,0.65,0.00}
\definecolor{orchid1}{rgb}{1.00,0.51,0.98}
\definecolor{orchid2}{rgb}{0.93,0.48,0.91}
\definecolor{orchid3}{rgb}{0.80,0.41,0.79}
\definecolor{orchid4}{rgb}{0.55,0.28,0.54}
\definecolor{orchid}{rgb}{0.85,0.44,0.84}
\definecolor{palegoldenrod}{rgb}{0.93,0.91,0.67}
\definecolor{palegreen}{rgb}{0.60,0.98,0.60}
\definecolor{paleturquoise}{rgb}{0.69,0.93,0.93}
\definecolor{paleviolet}{rgb}{0.86,0.44,0.58}
\definecolor{papayawhip}{rgb}{1.00,0.94,0.84}
\definecolor{peachpuff}{rgb}{1.00,0.85,0.73}
\definecolor{peru}{rgb}{0.80,0.52,0.25}
\definecolor{pink1}{rgb}{1.00,0.71,0.77}
\definecolor{pink2}{rgb}{0.93,0.66,0.72}
\definecolor{pink3}{rgb}{0.80,0.57,0.62}
\definecolor{pink4}{rgb}{0.55,0.39,0.42}
\definecolor{pink}{rgb}{1.00,0.75,0.80}
\definecolor{plum1}{rgb}{1.00,0.73,1.00}
\definecolor{plum2}{rgb}{0.93,0.68,0.93}
\definecolor{plum3}{rgb}{0.80,0.59,0.80}
\definecolor{plum4}{rgb}{0.55,0.40,0.55}
\definecolor{plum}{rgb}{0.87,0.63,0.87}
\definecolor{powderblue}{rgb}{0.69,0.88,0.90}
\definecolor{purple1}{rgb}{0.61,0.19,1.00}
\definecolor{purple2}{rgb}{0.57,0.17,0.93}
\definecolor{purple3}{rgb}{0.49,0.15,0.80}
\definecolor{purple4}{rgb}{0.33,0.10,0.55}
\definecolor{purple}{rgb}{0.63,0.13,0.94}
\definecolor{red1}{rgb}{1.00,0.00,0.00}
\definecolor{red2}{rgb}{0.93,0.00,0.00}
\definecolor{red3}{rgb}{0.80,0.00,0.00}
\definecolor{red4}{rgb}{0.55,0.00,0.00}
\definecolor{red}{rgb}{1.00,0.00,0.00}
\definecolor{rosybrown}{rgb}{0.74,0.56,0.56}
\definecolor{royalblue}{rgb}{0.25,0.41,0.88}
\definecolor{saddlebrown}{rgb}{0.55,0.27,0.07}
\definecolor{salmon1}{rgb}{1.00,0.55,0.41}
\definecolor{salmon2}{rgb}{0.93,0.51,0.38}
\definecolor{salmon3}{rgb}{0.80,0.44,0.33}
\definecolor{salmon4}{rgb}{0.55,0.30,0.22}
\definecolor{salmon}{rgb}{0.98,0.50,0.45}
\definecolor{sandybrown}{rgb}{0.96,0.64,0.38}
\definecolor{seagreen}{rgb}{0.18,0.55,0.34}
\definecolor{seashell1}{rgb}{1.00,0.96,0.93}
\definecolor{seashell2}{rgb}{0.93,0.90,0.87}
\definecolor{seashell3}{rgb}{0.80,0.77,0.75}
\definecolor{seashell4}{rgb}{0.55,0.53,0.51}
\definecolor{seashell}{rgb}{1.00,0.96,0.93}
\definecolor{sienna1}{rgb}{1.00,0.51,0.28}
\definecolor{sienna2}{rgb}{0.93,0.47,0.26}
\definecolor{sienna3}{rgb}{0.80,0.41,0.22}
\definecolor{sienna4}{rgb}{0.55,0.28,0.15}
\definecolor{sienna}{rgb}{0.63,0.32,0.18}
\definecolor{skyblue}{rgb}{0.53,0.81,0.92}
\definecolor{slateblue}{rgb}{0.42,0.35,0.80}
\definecolor{slategray}{rgb}{0.44,0.50,0.56}
\definecolor{slategrey}{rgb}{0.44,0.50,0.56}
\definecolor{snow1}{rgb}{1.00,0.98,0.98}
\definecolor{snow2}{rgb}{0.93,0.91,0.91}
\definecolor{snow3}{rgb}{0.80,0.79,0.79}
\definecolor{snow4}{rgb}{0.55,0.54,0.54}
\definecolor{snow}{rgb}{1.00,0.98,0.98}
\definecolor{springgreen}{rgb}{0.00,1.00,0.50}
\definecolor{steelblue}{rgb}{0.27,0.51,0.71}
\definecolor{tan1}{rgb}{1.00,0.65,0.31}
\definecolor{tan2}{rgb}{0.93,0.60,0.29}
\definecolor{tan3}{rgb}{0.80,0.52,0.25}
\definecolor{tan4}{rgb}{0.55,0.35,0.17}
\definecolor{tan}{rgb}{0.82,0.71,0.55}
\definecolor{thistle1}{rgb}{1.00,0.88,1.00}
\definecolor{thistle2}{rgb}{0.93,0.82,0.93}
\definecolor{thistle3}{rgb}{0.80,0.71,0.80}
\definecolor{thistle4}{rgb}{0.55,0.48,0.55}
\definecolor{thistle}{rgb}{0.85,0.75,0.85}
\definecolor{tomato1}{rgb}{1.00,0.39,0.28}
\definecolor{tomato2}{rgb}{0.93,0.36,0.26}
\definecolor{tomato3}{rgb}{0.80,0.31,0.22}
\definecolor{tomato4}{rgb}{0.55,0.21,0.15}
\definecolor{tomato}{rgb}{1.00,0.39,0.28}
\definecolor{turquoise1}{rgb}{0.00,0.96,1.00}
\definecolor{turquoise2}{rgb}{0.00,0.90,0.93}
\definecolor{turquoise3}{rgb}{0.00,0.77,0.80}
\definecolor{turquoise4}{rgb}{0.00,0.53,0.55}
\definecolor{turquoise}{rgb}{0.25,0.88,0.82}
\definecolor{violetred}{rgb}{0.82,0.13,0.56}
\definecolor{violet}{rgb}{0.93,0.51,0.93}
\definecolor{wheat1}{rgb}{1.00,0.91,0.73}
\definecolor{wheat2}{rgb}{0.93,0.85,0.68}
\definecolor{wheat3}{rgb}{0.80,0.73,0.59}
\definecolor{wheat4}{rgb}{0.55,0.49,0.40}
\definecolor{wheat}{rgb}{0.96,0.87,0.70}
\definecolor{whitesmoke}{rgb}{0.96,0.96,0.96}
\definecolor{white}{rgb}{1.00,1.00,1.00}
\definecolor{yellow1}{rgb}{1.00,1.00,0.00}
\definecolor{yellow2}{rgb}{0.93,0.93,0.00}
\definecolor{yellow3}{rgb}{0.80,0.80,0.00}
\definecolor{yellow4}{rgb}{0.55,0.55,0.00}
\definecolor{yellowgreen}{rgb}{0.60,0.80,0.20}
\definecolor{yellow}{rgb}{1.00,1.00,0.00}

\usepackage{hyperref}
\hypersetup{pdfpagemode=FullScreen,   
backref,
colorlinks=true,
citecolor=Bittersweet,
linkcolor=Bittersweet,
urlcolor=Bittersweet
}

\usepackage[labelsep=space]{caption}

\usepackage{epsfig}

\newtheorem{thm}{Theorem}
\newtheorem{mdef}{Definition}
\newtheorem{prop}{Proposition}

\begin{document}

\renewcommand{\eqref}[1]{(\ref{#1})}
\newcommand{\mb}[1]{\mathbf{#1}}
\newcommand{\mbb}[1]{\mathbb{#1}}
\newcommand{\R}{\mathbb{R}}
\newcommand{\mt}[1]{\mathrm{#1}}
\newcommand{\rv}{random variable}
\newcommand{\cqfd}{\hfill $\square$}

\title{Distance-based Depths for Directional Data}
\author{\large Giuseppe Pandolfo$^*$, Davy Paindaveine$^\dagger$ and Giovanni Porzio$^\ddagger$\\[3mm]
{\normalsize University of Naples Federico II$^*$,
 Universit\'{e} libre de Bruxelles$^\dagger$}, \\
{\normalsize and University of Cassino and Southern Lazio$^\ddagger$}
}
\date{}

\maketitle



\begin{abstract}
Directional data are constrained to lie on the unit sphere of~$\R^q$ for some~$q\geq 2$. To address the lack of a natural ordering for such data, depth functions have been defined on spheres. However, the depths available either lack flexibility or are so computationally expensive that they can only be used for very small dimensions~$q$. In this work, we improve on this by introducing a class of distance-based depths for directional data. Irrespective of the distance adopted, these depths can easily be computed in high dimensions too. We derive the main structural properties of the proposed depths and study how they depend on the distance used. We discuss the asymptotic and robustness properties of the corresponding deepest points. We show the practical relevance of the proposed depths in two applications, related to (i) spherical location estimation and (ii) supervised classification. For both problems, we show through simulation studies that distance-based depths have strong advantages over their competitors.
\end{abstract}


\section{Introduction}

Directional data analysis is relevant when the sample space is the unit hypersphere~$\mathcal{S}^{q-1}\linebreak :=\left\{x \in \mathbb{R}^{q}:x^{T}x = 1\right\}$ in $\mathbb{R}^{q}$, which occurs when observations are directions, axes, rotations, or cyclic events. Applications arise in numerous fields, including astronomy, earth sciences, biology, meteorology and political science; see \cite{GiHa2010} for an exemple in the latter field. Directional data analysis can also be exploited to study patterns of unit vectors in~$\mathbb{R}^{q}$, such as those encountered in text mining \citep{Hoetal2012}.

Statistically, analyzing and describing directional data requires tackling some interesting problems associated with the lack of a reference direction and with a sense of rotation not uniquely defined. Another important issue when dealing with such data is the lack of a natural ordering, which generates a special interest in depth functions on the sphere. Parallel to their role in the usual Euclidean case, directional depths are to measure the degree of centrality of a given spherical location with respect to a distribution on the sphere and to provide a center-outward ordering of spherical locations; see \cite{AgoRom2013a}.  

Depth concepts for directional data were first considered by \cite{Sma1987} and \cite{LiuSin1992}. Following the pioneering work of \cite{Sma1987}, \cite{LiuSin1992} popularized the concept of \textit{angular Tukey depth} (ATD), which is the directional analog of the celebrated \emph{halfspace depth} \citep{Tuk1975}. The same paper introduced two further depths for directional data, namely the \textit{angular simplicial depth} (ASD), which is the directional version of the \emph{simplicial depth} from \cite{Liu1990}, and the \emph{arc distance depth} (ADD), which is based on the concept of arc length distance. 

Unlike the ADD, the ATD and ASD have been studied and used in the literature. For instance, \cite{RouStr2004} investigated some of the properties of the ATD, while \cite{AgoRom2013a} considered some of the possible applications of the ASD and ATD. {\ttfamily{R}} packages are also available for these depths: the package {\ttfamily{depth}} (\citealp{Genetal2012}) allows to compute ATD values for~$q = 2$ or $3$, whereas the package {\ttfamily{localdepth}} (\citealp{AgoRom2013b}) implements specific functions for the evaluation of the ATD for~$q=2$, and of the ASD for an arbitrary~$q\geq 2$. 

The main drawback of both the ASD and ATD is the computational effort they require, especially for higher dimensions~$q$. The \emph{angular Mahalanobis depth} of \cite{Leyetal2014}, that is based on a concept of directional quantiles, is computationally much more affordable, but suffers from other disadvantages: it requires the preliminary choice of a spherical location functional and it is less flexible than the ASD/ATD in the sense that it produces rotationally symmetric depth contours, even if the underlying distribution is not rotationally symmetric.

On the one hand, depth functions for directional data are useful, yet on the other hand, they lack flexibility (and depend on some user's choice) or are computationally too demanding. In order to improve on this, this work introduces a new class of directional depth functions that is based on spherical distances and contains the ADD as a particular case. These depth functions are computationally feasible even in high dimensions and are generally more flexible. Distance-based directional depths show several other advantages over their ASD/ATD competitors: they take positive values everywhere on~$\mathcal{S}^{q-1}$ (but in the uninteresting case of a point mass distribution), whereas the ASD/ATD can take zero values (which is undesirable when performing supervised classification). Further advantages of the proposed distance-based depths is that they typically do not provide ties in the sample case (whereas ties are unavoidable for the ASD/ATD, due to their step function nature) and that they do not require any assumption on the underlying distribution (unlike the angular Mahalanobis depth that, when based on the spherical mean, is not defined for zero-mean distributions).

The paper is organized as follows. In Section~\ref{secdef}, we introduce the proposed class of distance-based depth functions for directional data, and we consider three particular cases, namely the arc distance depth (ADD), the cosine distance depth (CDD) and the chord distance depth (ChDD). In Section~\ref{sec:StructProperties}, we derive the main structural properties of the proposed depths and study how they depend on the distance used. In Section~\ref{secillu}, we compare the various depths considered for several empirical distributions on the circle ($q=2$), which also allows us to illustrate the theoretical results of Section~\ref{sec:StructProperties}. In Section~\ref{sec:DistrProperties}, we discuss the asymptotic and robustness properties of the proposed concepts. In Section~\ref{secSimu}, we show the practical relevance of the distance-based depths in two applications, related to (i) spherical location estimation (Section~\ref{secSimusub1}) and (ii) supervised classification (Section~\ref{secSimusub2}). For both problems, we perform simulations that show the advantages of the proposed depths over their competitors. Final comments are provided in Section~\ref{secfinal}. Finally, an appendix collects technical proofs.


\section{Distance-based depths for directional data}
\label{secdef}

\noindent In Definition~\ref{defclass} below, we introduce a class of depths on the unit sphere~$\mathcal{S}^{q-1}$. A particular member of this class will be obtained by fixing a particular (bounded) distance~$d(\cdot,\cdot)$ on~$\mathcal{S}^{q-1}$. For such a distance, $d^{\rm sup}:=\sup \{ d(\theta, \psi): \theta,\psi\in \mathcal{S}^{q-1}\}$ will throughout denote the upper bound of the distance between any two points on~$\mathcal{S}^{q-1}$.   

\vspace{-.3cm}
\begin{mdef}[Directional distance-based depths] 
\label{defclass}
Let~$d(\cdot,\cdot)$ be a bounded distance on $\mathcal{S}^{q-1}$ and $H$ be a distribution on~$\mathcal{S}^{q-1}$. Then the 
\emph{directional $d$-depth of~$\theta(\in \mathcal{S}^{q-1})$ with respect to~$H$} is
\begin{align}
\label{eq:class}
D_{d}\left(\theta, H\right) 
:= 
d^{\rm sup} 
- 
E_H[d(\theta, W)]
,	
\end{align}
where~$E_H$ is the expectation under the assumption that~$W$ has distribution~$H$.
\end{mdef}
\vspace{-.3cm}
While, in principle, any distance~$d$ can be used in this definition, it is natural to consider distances that are \emph{rotation-invariant} in the sense that~$d(O\theta,O\psi)=d(\theta,\psi)$ for any~$\theta,\psi\in\mathcal{S}^{q-1}$ and any $q\times q$ orthogonal matrix~$O$. As we show for the sake of completeness in the appendix (see Proposition~\ref{prodistinv}), any rotation-invariant distance~$d$ is of the form
$$
d(\theta,\psi)=d_{\delta}(\theta,\psi)=\delta(\theta'\psi)
$$ 
for some function~$\delta:[-1,1]\to\R^+$. The standard distance axioms impose that~$\delta(1)=0$ but do not impose that~$\delta$ is monotone non-increasing (unexpectedly, the triangle inequality may hold without this monotonicity condition). All classical choices, however, are monotone non-increasing; these include the \emph{arc length distance}~$d_{\rm arc}$ and the \emph{cosine distance}~$d_{\cos}$, that are associated with~$\delta(t)=\delta_{\rm arc}(t)=\arccos t$ and~$\delta(t)=\delta_{\cos}(t)=1-t$, respectively. Another rotation-invariant distance for which this monotonicity condition holds is the \emph{chord distance}~$d_{\rm chord}$ defined through~$d_{\rm chord}(\theta,\psi)=\|\theta-\psi\|=\sqrt{2(1-\theta'\psi)}
\linebreak =:\delta_{\rm chord}(\theta'\psi)$. Throughout, we will denote the corresponding \emph{arc distance depth} (ADD), \emph{cosine distance depth} (CDD) and \emph{chord distance depth} (ChDD) as~$D_{\rm arc}$,~$D_{\cos}$ and~$D_{\rm chord}$, respectively.  

The ADD is the arc distance depth introduced by \cite{LiuSin1992}. For the CDD, a direct computation yields
\begin{equation}
	\label{cosinexplic}
D_{\cos}(\theta, H)
= 
2
- 
E_{H}[1-\theta' W]
= 
1
+ 
\theta' E_{H}[W]
.
\end{equation}
Under the assumption that~$E_{H}[W]$ is non-zero, this rewrites
$D_{\cos}(\theta, H) 
= 
1
+ 
\|E_{H}[W]\| 
\linebreak (\theta' \mu_H)
$,
where~$\mu_H:=E_{H}[W]/\|E_{H}[W]\|$ is the spherical mean of~$H$. This shows that the CDD is then in a one-to-one relationship with the \emph{angular Mahalanobis depth} of \cite{Leyetal2014}, provided that the location functional needed in the latter is chosen as the spherical mean. We stress, however, that, unlike the angular Mahalanobis depth, the CDD does not require choosing a location functional on the sphere and is defined also in cases where~$\mu_H=0$. To the best of our knowledge, the ChDD has not been considered in the literature. 


\section{Structural properties}
\label{sec:StructProperties}

In this section, we derive the main properties of a generic directional $d$-depth. We start with the following  invariance result. 

\begin{thm}{\textbf{(Rotational invariance)}} 
\label{thmrotainv}
Let~$d=d_\delta$ be a rotation-invariant distance and~$H$ be a distribution on~$\mathcal{S}^{q-1}$. Then $D_{d_{\delta}}(\theta, H)$ is a rotation-invariant depth, in the sense that~$D_{d_{\delta}}(O\theta, H_O)=D_{d_{\delta}}(\theta, H)$ for any $q\times q$ orthogonal matrix~$O$, where~$H_O$ denotes the image of~$H$ by the transformation~$x\mapsto Ox$, that is, $H_O$ is the distribution of~$OW$ when~$W$ has distribution~$H$.
\end{thm}

A corollary is that if~$H$ is rotationally symmetric about~$\theta_0$ in the sense that~$H_O=H$ for any $q\times q$ orthogonal matrix~$O$ fixing~$\theta_0$, then $d_\delta (O\theta,H)=d_\delta (\theta, H)$ for any such $O$. In particular, for any~$\alpha$, the $\alpha$-depth region --- that, as usual, is defined as the collection of~$\theta$ values with a depth larger than or equal to~$\alpha$ --- is invariant under rotations fixing~$\theta_0$, hence reflects the symmetry of the distribution~$H$ about~$\theta_0$. 

In contrast, parallel to the angular Mahalanobis depth of \cite{Leyetal2014}, the CDD provides symmetric depth regions of this form for any $H$, i.e, irrespectively of the fact that~$H$ is rotationally symmetric or not. This follows from the comments at the end of Section~\ref{secdef}.

\begin{thm}{\textbf{(Continuity)}} 
\label{thcontinuity}
Assume that the distance~$d$ is continuous; if~$d=d_\delta$, then this is equivalent to assuming that~$\delta:[-1,1]\to \R^+$ is continuous. Let~$H$ be a distribution on~$\mathcal{S}^{q-1}$. Then, (i) the mapping~$\theta\mapsto D_{d}(\theta, H)$ is continuous on~$\mathcal{S}^{q-1}$; (ii) there exists~$\theta_{d}(H)\in\mathcal{S}^{q-1}$ such that $D_{d}(\theta_d(H),H)=\sup_{\theta\in\mathcal{S}^{q-1}} D_{d}(\theta,H)$. 
\end{thm}

Note that the continuity result in Theorem~\ref{thcontinuity}(i) holds without any assumption on~$H$, hence will also hold in the empirical case. Theorem~\ref{thcontinuity}(ii) guarantees the existence of a $D_{d}$-deepest point~$\theta_{d}(H)$. The deepest point (or collection of deepest points) typically depends on the distance~$d$ adopted. For the CDD, the deepest point is the spherical mean, provided that~${\rm E}_H[W]\neq 0$, whereas the deepest point for the ADD is the spherical median of \cite{Fis1985}, which reduces to the circular median (\citealp{MarJup2000},  p.~30) in dimension~$q=2$. 
This is in line with the Euclidean case where deepest points typically depend on the depth considered and may be multivariate medians (e.g., Tukey's halfspace or Liu's simplicial deepest points) or mean vectors (e.g., the zonoid of \cite{KosMos1997} or the moment-based Mahalanobis deepest points). 

The deepest point may not be unique; for the uniform distribution on $\mathcal{S}^{q-1}$, for instance, any rotation-invariant distance-based depth will be constant over the sphere (this readily follows from Theorem~\ref{thmrotainv}). 
This lack of unicity also holds in the Euclidean case, where the barycentre of the collection $\mathcal{C}$ of deepest points is often taken as its unique representative; for most depths, it then follows from the convexity of the depth regions (which guarantees convexity of $\mathcal{C}$) that this barycentre indeed has maximal depth.  
It is interesting to note that directional depths 
are fundamentally different in this respect, as no such convexity arguments can be used. The particular nature of the sample space may induce depth regions that are even disconnected. This may occur for some multimodal distributions~$H$; an example is given in Section~\ref{secillu}. 
In contrast, note that, for~$D_{\cos}$, the collection of deepest points is either $\left\{\mu_{H}\right\}$, when ${\rm E}_H[W]\neq 0$ , or $\mathcal{S}^{q-1}$, when ${\rm E}_H[W]=0$, and hence it is always spherically convex.

It is desirable that if the distribution~$H$ on~$\mathcal{S}^{q-1}$ has an ``indisputable" location centre~$\theta_0$, then the deepest point~$\theta_d(H)$ is unique and coincides with~$\theta_0$. The following theorem provides such a Fisher consistency result. 

\begin{thm}{\textbf{(Fisher consistency under monotone rotational symmetry)}} 
\label{thFishconsist}
Assume that the rotation-invariant distance~$d=d_\delta$ is based on a monotone strictly decreasing function~$\delta:[-1,1]\to \R^+$. Assume that the distribution~$H$ on~$\mathcal{S}^{q-1}$ admits a density of the form~$x\mapsto c_{q,h} h(x'\theta_0)$ for some~$\theta_0\in\mathcal{S}^{q-1}$ and some monotone strictly increasing function~$h:[-1,1]\linebreak \to \R^+$. Then, $\theta\mapsto D_{d_{\delta}}(\theta, H)$ is a monotone strictly increasing function of~$\theta'\theta_0$, so that $\theta\mapsto D_{d_{\delta}}(\theta, H)$ is uniquely maximized at~$\theta_0$. 
\end{thm}

Theorem~\ref{thFishconsist} ensures that the ADD-, CDD-, and ChDD-deepest points are equal and coincide with the modal location $\theta_0$ of $H$ in case the latter admits a density of the form given in the theorem.
The monotonicity result entails that, irrespective of the distance~$d_\delta$ used, the depth regions are of the form~$\{\theta\in\mathcal{S}^{q-1}: \theta'\theta_0\geq c\}$. 

In this setup, the maximal depth,~$\max_{\theta\in\mathcal{S}^{q-1}} D_{d_\delta}(\theta,H)$, measures the concentration of~$H$, as showed in the following theorem. 

\begin{thm}{\textbf{(Maximal depth as a concentration measure)}} 
\label{propconcentr}
Assume that the rotation-invariant distance~$d=d_\delta$ is based on a monotone strictly decreasing function~$\delta:[-1,1]\to \R^+$. Assume that the distribution~$H_\kappa$ on~$\mathcal{S}^{q-1}$ admits the density~$x\mapsto c_{q,\kappa,h} h(\kappa x'\theta_0)$ for some~$\theta_0\in\mathcal{S}^{q-1}$ and some monotone strictly increasing and differentiable function~$h:\R\to \R^+$ such that
$
t \mapsto t\,\frac{d}{dt}\log h(t)
$
is monotone strictly increasing. Then the maximal depth~$D_{d_{\delta}}(\theta_0, H_\kappa)$ is a strictly increasing function of~$\kappa$.  
\end{thm}
In Theorem \ref{propconcentr}, $\kappa$ plays the role of a concentration parameter; typically, the larger~$\kappa$,  the more concentrated the probability mass is about the modal location~$\theta_0$. Since the maximal depth is a strictly increasing function of~$\kappa$, it is itself a concentration (or spread) measure. Note that the assumption that $t \mapsto t\,\frac{d}{dt}\log h(t)
$ is monotone strictly increasing holds in particular if~$h$ is $\log$-convex, so that the result applies for von Mises--Fisher (vMF) distributions that are obtained for~$h\left(u\right)=\exp\left(u\right)$. While Theorem~\ref{propconcentr} restricts to rotationally symmetric distributions, the maximal cosine distance depth~$\max_{\theta\in\mathcal{S}^{q-1}} D_{\cos}(\theta,H)=1+\|E_{H}[W]\|$ is, irrespective of~$H$, related to the ``spherical variance" (\citealp{MarJup2000}, p.~164), that is, to the mean resultant length~$\|E_H[W]\|$ of~$W$.

We conclude this section by stating a property showing that the proposed depths may inherit anti-symmetry properties of the distances on which they are based. More precisely, we have the following result which is restricted to rotationally-invariant distances, although a similar result can be stated for an arbitrary distance~$d$. 
\setlength{\topsep}{1.5em}
\begin{thm}{\textbf{(Anti-symmetry)}} 
\label{thskewsym}
Assume that the rotation-invariant distance~$d=d_\delta$ is based on a function~$\delta:[-1,1]\to \R^+$ that is anti-symmetric about~$0$, i.e., $\delta(-t)+\delta(t)=\delta(-1)$. Let~$H$ be a distribution on~$\mathcal{S}^{q-1}$. Then, 
\begin{enumerate}[label=(\roman*), topsep=0pt, itemsep=-1ex]
\item $\theta\mapsto D_{d_{\delta}}(\theta, H)$ is anti-symmetric on~$\mathcal{S}^{q-1}$ in the sense that
\begin{equation*}
D_{d_{\delta}}(-\theta, H)=d_{\delta}^{\sup}-D_{d_{\delta}}(\theta,H);
\end{equation*}
\item If~$\theta_0$ has maximal depth, then~$-\theta_0$ has minimal depth.  
\end{enumerate}
\end{thm}
The arc length and cosine distances are based on anti-symmetric functions~$\delta$, but the chord distance is not. 
If~$\delta$ is anti-symmetric, then an antipodally symmetric distribution~$H \in \mathcal{S}^{q-1}$, under which~$H(-B) = H(B)$ for any measurable set~$B$ on~$\mathcal{S}^{q-1}$~, leads to a depth function $\theta\mapsto D_{d_{\delta}}(\theta,H)$ that is constant. This is another property contrasting sharply with the Euclidean case, where no distribution will provide a constant depth function. 
To show why the claim on the constancy holds true, consider an arbitrary measurable set~$B\subset\mathcal{S}^{q-1}$ such that~$\mathcal{S}^{q-1}=(-B)\cup B$ and~$(-B)\cap B=\emptyset$. Then, using the antipodal symmetry of~$H$ and the antisymmetry of~$\delta$, we obtain 
\begin{eqnarray*}
 D_{d_{\delta}}(\theta,H)
&=&
 \delta(-1)
-
\int_{-B}
\delta(\theta'w)
\,
dH(w)
-
\int_{B}
\delta(\theta'w)
\,
dH(w)
\\[2mm]
&= & 
 \delta(-1)
-
\int_{B}
\delta(-\theta'w)
\,
dH(w)
-
\int_{B}
\delta(\theta'w)
\,
dH(w)
\\[2mm]
&= & 
 \delta(-1)
-
\int_{B}
\delta(-1)
\,
dH(w)
=
\frac{ \delta(-1)}{2}
\cdot
\end{eqnarray*}
An interesting question is whether or not antipodal symmetry of~$H$ is also a necessary condition for the constancy of~$\theta\mapsto D_{d_{\delta}}(\theta,H)$ with an anti-symmetric function~$\delta$. While \cite{LiuSin1992} proved that this is indeed the case for the ADD in dimension~$q=2$ under the assumption that~$H$ admits a density, it is not the case for any~$\delta$ function. For instance, for the CDD, it directly follows from~(\ref{cosinexplic}) that~$\theta\mapsto D_{\cos}(\theta, H)$ is constant if and only if $E_{H}[W]=0$, which shows that antipodal symmetry is not a necessary condition for the constancy of~$D_{\cos}$. 


\section{Illustrations}
\label{secillu}

This short section illustrates the theoretical results of the previous section for three empirical distributions on the circle~$\mathcal{S}^1$; we restrict to the circle to allow for a visual comparison of the various depths. Denoting as $H^{\rm vMF}_{\alpha,\kappa}$ the vMF distribution on~$\mathcal{S}^1$ with modal location~$\theta=(\cos\alpha,\sin\alpha)'$ and concentration~$\kappa$, the three empirical distributions considered are associated with a random sample of size~$n=500$ from each of the following distributions:
$H_1=H^{\rm vMF}_{\pi,2}$ (unimodal case),
$H_2=\frac{1}{2} H^{\rm vMF}_{\frac{3\pi}{4},5}+\frac{1}{2} H^{\rm vMF}_{\frac{5\pi}{4},5}$ (bimodal symmetric case), 
$H_3=\frac{1}{2} H^{\rm vMF}_{\frac{5\pi}{9},7}+\frac{1}{2} H^{\rm vMF}_{\frac{13\pi}{9},17}$ (bimodal asymmetric case).

\begin{figure}[h!]
\captionsetup{font=scriptsize}
\begin{center}
\includegraphics[width=1.00\textwidth]{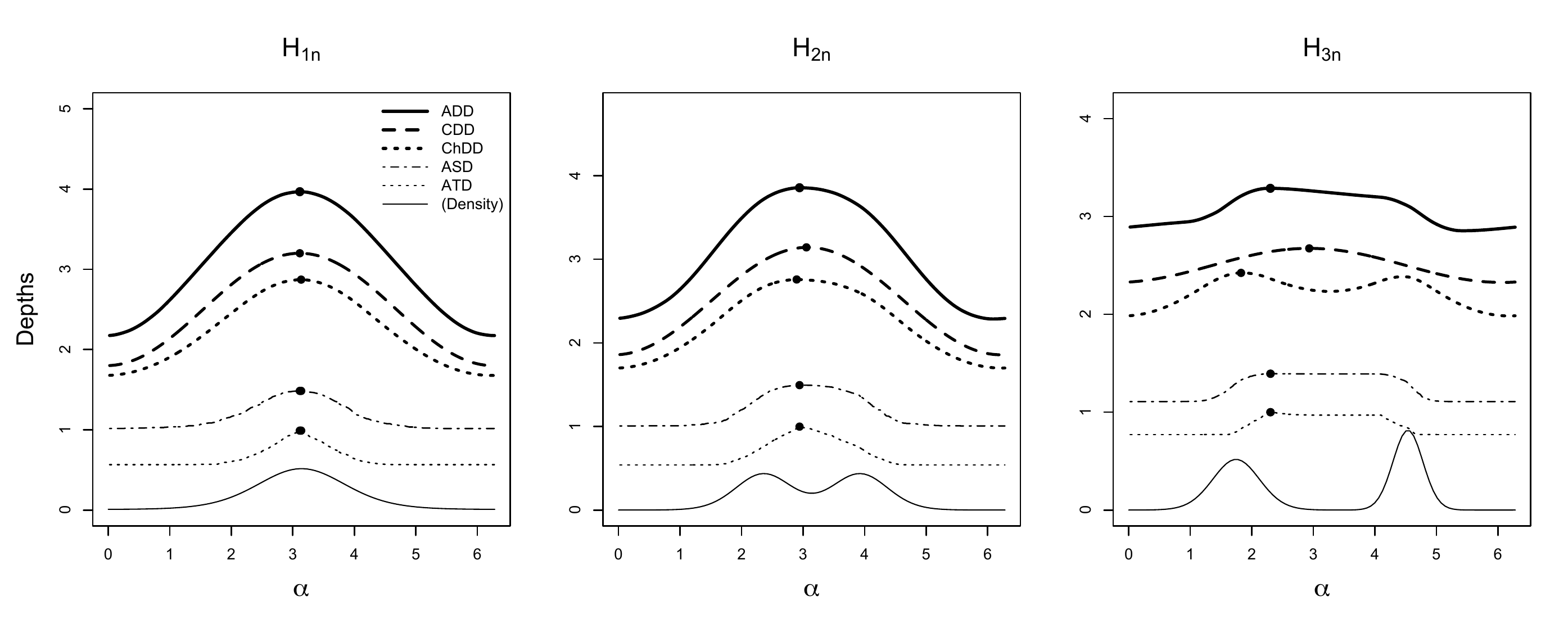}
\end{center}
\vspace{-5mm}
\caption{Plots of the depth mapping~$\alpha\mapsto D({\cos \alpha \choose \sin  \alpha},H_{\ell n})$, for the distance-based depths ADD, CDD and ChDD,
as well as the angular simplicial depth (ASD) and angular Tukey depth (ATD),
 and the empirical distributions~$H_{\ell n}$, $\ell=1,2,3$ described in Section~\ref{secillu} (for easier visualization, depth values were actually multiplied by~1.5 for distance-based depths, by~1 for the ASD, and by~0.5 for the ATD). Deepest points are maked by a black dot. The parent density is also plotted in each case. 
}
\label{Figillu}
\end{figure}

For each of the resulting empirical distributions~$H_{\ell n}$, $\ell=1,2,3$,   Figure~\ref{Figillu} provides plots of the distance-based depths ADD, CDD and ChDD, as well as the competing angular simplicial depth (ASD) and angular Tukey depth (ATD). The ASD and ATD were computed through the packages {\ttfamily{localdepth}} and {\ttfamily{depth}}, respectively. The distance-based depths were computed by means of {\ttfamily{R}} functions written by the authors. Simulated data and their graphical representations were obtained through the {\ttfamily{R}} package {\ttfamily{circular}} (\citealp{LunAgo2013}), which is a standard reference to work with data on the unit circle. 

For~$H_{1n}$, all distance-based depth functions are monotonically strictly decreasing from their deepest point ($\approx \pi$) and do so in a symmetric way, which is in accordance with Theorems~\ref{thmrotainv} and~\ref{thFishconsist}. These depths functions are also continuous; see Theorem~\ref{thcontinuity}. In contrast, the ATD is constant outside the interval of length~$\pi$ centered at its deepest point, which holds for any distribution on the circle (\citealp[Proposition~4.6.]{LiuSin1992}), and both the ASD and ATD are piecewise constant functions. The center-outward rankings provided by the ASD and ATD therefore yield many ties and are more rough than those given by distance-based depths.  
For the symmetric bimodal distribution~$H_{2n}$, all depth functions are unimodal, hence fail to capture the bimodality of the distribution, which is not a problem since depths are not density measures but rather centrality measures. 
In contrast with the Euclidean case, some directional depths may exhibit multimodality, as it is the case for the ChDD for the distribution $H_{3n}$, where modes are more separated than in $H_{2n}$; (\ref{cosinexplic}) entails that the CDD will never exhibit such a multimodal pattern. In this last example, the depth functions reflect the asymmetry of the distribution and do not identify the same deepest point; in particular, the CDD is maximized at the spherical mean, whereas the ADD is maximized at the circular median (\citealp{MarJup2000}, p.~20), and so are the ASD and ATD.


\section{Asymptotic and robustness properties}
\label{sec:DistrProperties}

In this section, we present asymptotic results for the distance-based depths introduced in Definition~\ref{defclass} and for the corresponding deepest points, as well as a robustness result regarding the breakdown point of these. We start with a Glivenko-Cantelli-type result. 

\begin{thm}{\textbf{(Uniform almost sure consistency)}} 
\label{thunifconsistency}
Let~$d$ be a bounded and continuous distance on~$\mathcal{S}^{q-1}$ and~$H$ be a distribution on~$\mathcal{S}^{q-1}$. Denote as~$H_n$ the empirical distribution associated with a random sample of size~$n$ from~$H$. Then 
$$
\sup_{\theta\in\mathcal{S}^{q-1}} 
\big|D_{d}(\theta, H_n)-D_{d}(\theta, H)\big|
\to 
0
$$
almost surely as~$n\to\infty$. 
\end{thm}

This result implies that we may explore empirically the properties of~$D_d(\theta,H)$ by considering the corresponding sample depth function~$D_d(\theta,H_n)$ for a large~$n$. This justifies a posteriori the illustration of Theorem~\ref{thFishconsist} in the previous section.
The following asymptotic normality is a direct result of the central limit theorem.

\begin{thm}{\textbf{(Asymptotic normality of sample depth)}} 
\label{thasnormD}
Let~$d$ be a bounded distance on~$\mathcal{S}^{q-1}$ and~$H$ be a distribution on~$\mathcal{S}^{q-1}$. Denote as~$H_n$ the empirical distribution associated with a random sample of size~$n$ from~$H$. Then as~$n\to\infty$,
$
\sqrt{n}(D_{d}(\theta, H_n)-D_{d}(\theta, H))
$
converges weakly to the normal distribution with mean zero and variance~${\rm Var}_H[d(\theta,W)]$. 
\end{thm}

We turn to asymptotic and robustness results for deepest points. The following strong consistency result requires that the deepest point is uniquely defined, as it is in Theorem~\ref{thFishconsist}. 

\begin{thm}{\textbf{(Almost sure consistency of the deepest point)}} 
\label{thdeepestconsistency}
Let~$d$ be a bounded and continuous distance on~$\mathcal{S}^{q-1}$ and~$H$ be a distribution on~$\mathcal{S}^{q-1}$. Assume that the deepest point~$\theta_d(H)$ is unique. Denote as~$H_n$ the empirical distribution associated with a random sample of size~$n$ from~$H$, and let~$\theta_d(H_n)$ be an arbitrary deepest point with respect to~$H_n$. Then
$$
\theta_d(H_n)
\to
\theta_d(H)
$$
almost surely as~$n\to\infty$.  
\end{thm}

Constructing confidence zones for~$\theta_d(H)$ requires the availability of the asymptotic distribution of~$\theta_d(H_n)$. Since~$\theta_d(H_n)$ is an $M$-estimator for a location parameter on~$\mathcal{S}^{q-1}$, its asymptotic distribution can easily be obtained from the results of \cite{KoCha1993}, at least under rotationally symmetric distributions. 
We do not pursue in this direction here.

Since deepest points are commonly used as robust location estimators, it is natural to investigate their robustness, and we therefore end this section by deriving a result on their breakdown point (BDP). In the directional setup considered, the classical BDP concept (\cite{Hametal1986}, pp.~97-98) is not suitable, and we adopt the directional concept of \cite{LiuSin1992} defining the BDP of the (more generally, of a) deepest point~$\theta_d(H)$ as the infimum 
of~$\varepsilon$ such that, for some contaminating distribution $G$ on~$\mathcal{S}^{q-1}$, $-\theta_d(H)$ is a deepest point of~$D_d(\theta, H_{\epsilon})$ with~$H_{\epsilon}:= \left(1 - \epsilon \right) H + \epsilon G$. The following result extends to an arbitrary distance~$d$ the lower bound result obtained in \cite{LiuSin1992} for the arc length distance.

\begin{thm}{\textbf{(Breakdown point of deepest points)}}
\label{thmbdp}
Let $d$ be a bounded distance on~$\mathcal{S}^{q-1}$ and~$H$ be a distribution on~$\mathcal{S}^{q-1}$. Let~$\theta_d(H)$ be a deepest point of~$D_{d}(\theta,H)$. Then the breakdown point of~$\theta_d(H)$ is larger than or equal to~$(D_{d}(\theta_d(H),H) - D_{d}(-\theta_d(H),H))\linebreak/(2d^{\rm sup})$.
\end{thm}

To investigate how the distance~$d$ affects the lower bound, we consider the important case of vMF distributions. If~$H^{\rm vMF}_{q,\theta_0,\kappa}$ denotes the vMF($\theta_0,\kappa$) distribution on~$\mathcal{S}^{q-1}$, then, for a rotation-invariant distance~$d_\delta$ that is decreasing in the sense of Theorem~\ref{thFishconsist}, we have~$\theta_{d_\delta}(H^{\rm vMF}_{q,\theta_0,\kappa})=\theta_0$ and  
$$
D_{d_\delta}(\pm\theta_0,H^{\rm vMF}_{q,\theta_0,\kappa})
=
d^{\rm sup}_\delta
-
\frac{ 
\int_{-1}^1
\delta(\pm v)
(1-v^2)^{(q-3)/2}
\exp(\kappa v) 
\,
dv
}
{\int_{-1}^1
(1-v^2)^{(q-3)/2}
\exp(\kappa v) 
\,
dv
} 
,
$$
which allows us to evaluate the lower bound from Theorem~\ref{thmbdp}. 

\begin{figure}[h!]
\captionsetup{font=scriptsize}
\begin{center} 
\makebox[\textwidth]{%
\includegraphics[width=\textwidth]{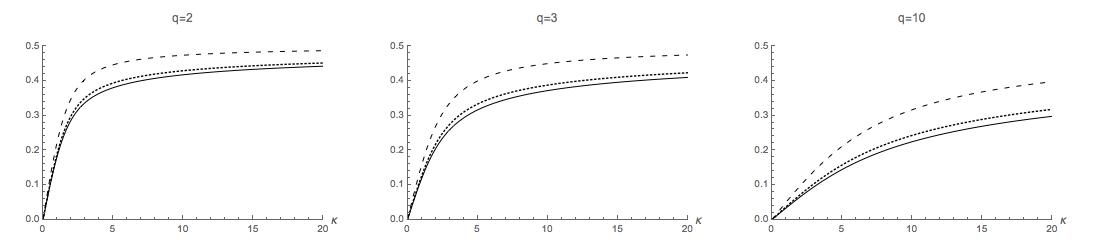}
}
\end{center}
\vspace{-1em}
\caption{Plots of the lower bound in Theorem~\ref{thmbdp}, for various dimensions~$q$ and for the ADD $(\rule[0.5ex]{0.4cm}{0.2pt})$, CDD $(\scriptsize\textendash\ \textendash)$, and ChDD $(\cdots)$, as a function of the concentration~$\kappa$ of the underlying vMF distribution on~$\mathcal{S}^{q-1}$.}
\label{BDPfig}
\end{figure}

Figure~\ref{BDPfig} plots this lower bound as a function of~$\kappa$ for various dimensions~$q$ and for the ADD, CDD and ChDD. Clearly, irrespective of the dimension and the distance, the lower bound is arbitrarily small for arbitrarily small values of~$\kappa$ and goes to 50\% as~$\kappa$ goes to infinity. The lower bound decreases as the dimension~$q$ increases. More importantly, for vMF distributions, the CDD-deepest point, namely the spherical mean,
provides a larger lower bound than the ADD- and CHDD-deepest ones do.


\section{Applications}
\label{secSimu}

We present two applications, which are related to spherical location estimation and supervised classification. 


\subsection{Spherical location estimation}
\label{secSimusub1}

Depth functions find applications in robust statistics, with the deepest point considered as a robust location estimator.

\begin{figure}[h!]
\captionsetup{font=scriptsize} 
\begin{center}  
\makebox[\textwidth]{%
\includegraphics[width=\textwidth]{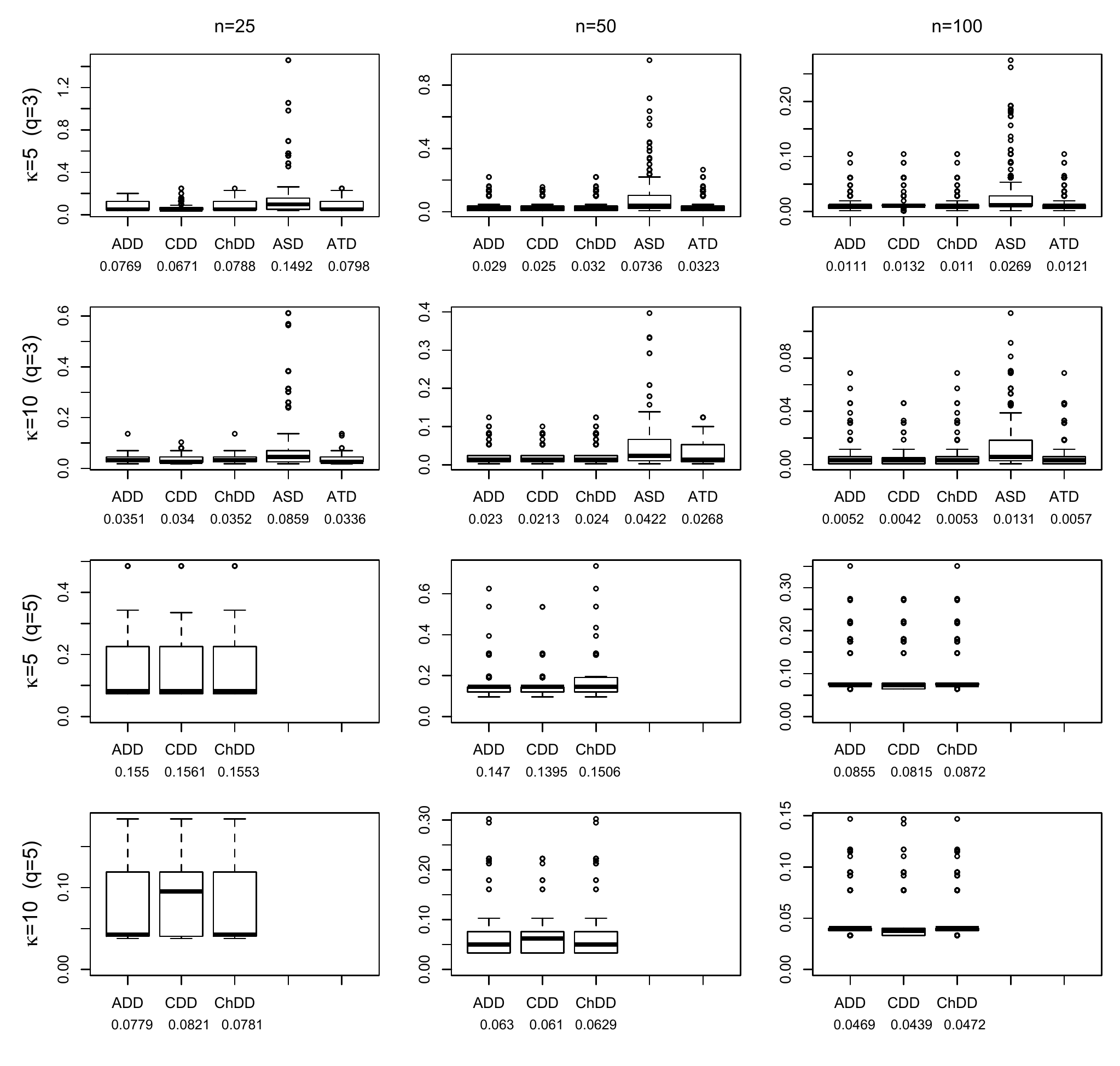}
}
\end{center}
\vspace{-5mm}
\caption{Boxplots, for~$q\in\{3,5\}$, $n\in\{25,50,100\}$ and~$\kappa\in\{5,10\}$, of the squared errors~${\rm SE}_{m}$, $m=1,\ldots,M$ (see~(\ref{SEdef})) of various depth-based estimators of~$\theta$ obtained from $M=500$ independent random samples of size~$n$ from the vMF distribution~$H^{\rm vMF}_{q,\theta,\kappa}$ with location~$\theta=e_q$ (the last vector of the canonical basis of~$\R^q$). The estimators considered are the ADD-, CDD- and ChDD-deepest points, as well as (due to computational issues, for dimension~$q=3$ only) the deepest points associated with the ASD and ATD. In each case, the corresponding mean square error~${\rm MSE} = (1/M) \sum_{m=1}^M {\rm SE}_{m}$ is provided.} 
\label{Efficiencyfig}
\end{figure}

For this reason, we conducted a simulation study to investigate the efficiency and robustness properties of the deepest points associated with the proposed distance-based depths, and to compare them with those of the competing ASD- and ATD-deepest points
We start with efficiency properties. For any combination of a dimension~$q\in\{3,5\}$, a sample size~$n\in\{25,50,100\}$ and a concentration~$\kappa\in\{5,10\}$, we generated $M=500$ independent random samples of size~$n$ from the distribution~$H^{\rm vMF}_{q,\theta,\kappa}$, where~$\theta=e_q$ is the last vector of the canonical basis of~$\R^q$. 
For each estimator~$\hat{\theta}$ of~$\theta$ considered, this leads to estimates~$\hat{\theta}_{1},\ldots,\hat{\theta}_{M}$. 
Figure~\ref{Efficiencyfig} provides boxplots of the resulting squared errors 
\begin{equation}
{\rm SE}_{m} 
= 
 \|\hat{\theta}_{m}- \theta\|^2
= 
2 (1 - \hat{\theta}_{m}' \theta)
,
\qquad
m=1,\ldots, M,
\label{SEdef}	
\end{equation}
and indicates the resulting mean square errors~${\rm MSE} = (1/M) \sum_{m=1}^M {\rm SE}_{m}$. The computational burden for the ASD- and ATD-deepest points is so prohibitive that these were considered for dimension~$q=3$ only.

Results indicate that, in dimension~$q=3$, the estimators associated with distance-based depths slightly dominate their ATD competitor and outperform their ASD one. As expected, the CDD-deepest point, that is the maximum likelihood estimator in the distributional setup considered, is in most cases the most efficient estimator. 
In dimension~$q=5$, where the ASD/ATD estimators could not be computed, the distance-based depths perform similarly. On the other hand, while the CDD estimator slightly dominates at all sample sizes in dimension $q = 3$, it dominates only at the largest considered sample size in dimension $q=5$.

\begin{figure}[h!]
\captionsetup{font=scriptsize}
\begin{center}  
\makebox[\textwidth]{%
\includegraphics[width=\textwidth]{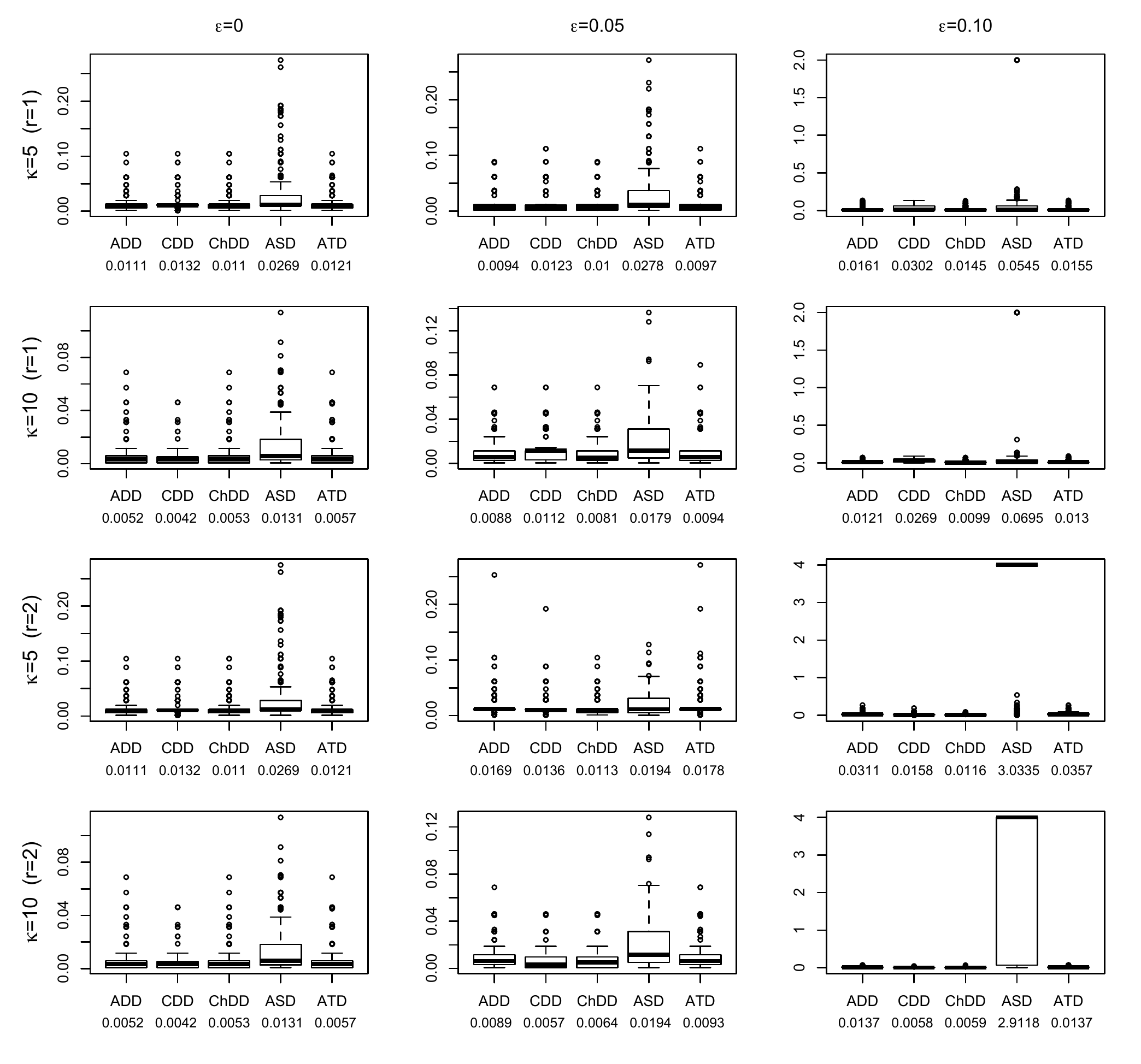}
}
\end{center}
\vspace{-5mm}
\caption{Boxplots, for~$q=3$, $\varepsilon\in\{0,0.05,0.10\}$ and~$\kappa\in\{5,10\}$, of the squared errors~${\rm SE}_{m}$, $m=1,\ldots,M$ (see~(\ref{SEdef})) of various depth-based estimators of~$\theta$ obtained from $M=500$ independent random samples of size~$n=100$ from the contaminated distribution~$(1-\varepsilon)H^{\rm vMF}_{q,\theta,\kappa}+\varepsilon \Delta_{\theta_r}$, where~$\theta$ is the last vector of the canonical basis of~$\R^q$, $\Delta_{\psi}$ denotes the point mass distribution at~$\psi$, and where $\theta_1$ (resp., $\theta_2$) is an orthogonal point to~$\theta$: $r=1$ (resp., the antipodal point to~$\theta$: $r=2$). The estimators considered are the ADD-, CDD- and ChDD-deepest points, as well as the deepest points associated with the ASD and ATD. In each case, the corresponding mean square error~${\rm MSE} = (1/M) \sum_{m=1}^M {\rm SE}_{m}$ is provided.} 
\label{Robustnessfig}
\end{figure}

We now turn to the investigation of robustness properties for which we restricted to dimension~$q=3$. For any combination of a contamination level~$\varepsilon\in\{0,0.05,0.10\}$ and a concentration~$\kappa\in\{5,10\}$, we generated $M=500$ independent random samples of size~$n=100$ from the contaminated distributions~$(1-\varepsilon)H^{\rm vMF}_{q,\theta,\kappa}+\varepsilon \Delta_{\theta_r}$, $r=1,2$, where~$\theta$ is set as~$e_q$, $\theta_1=e_{q-1}$, $\theta_2=-\theta$, $\Delta_{\psi}$ denotes the point mass distribution at~$\psi$. Hence, $r=1,2$ refers to contamination at an orthogonal point to~$\theta$ and at the antipodal point to~$\theta$, respectively. In each sample, the deepest points of the same five depths as in Figure~\ref{Efficiencyfig} were computed. The resulting boxplots of squared errors~${\rm SE}_m$ for $m=1,\ldots,M$ and the mean squared errors (${\rm MSE}$) are provided in Figure~\ref{Robustnessfig}. 

The results show that the estimators associated with distance-based depths enjoy good robustness properties. In particular, irrespective of the contamination level~$\varepsilon$ and the type of contamination, the ADD, CDD and ChDD estimators outperform the ASD one in terms of robustness. The domination over the ATD estimator is less significant.


\subsection{Supervised classification}
\label{secSimusub2}

Classification has been one of the most successful applications of statistical depth in the last decade, both for multivariate and functional data. While some proposals were based  on the use of local depth concepts (\citealp{PaiVanB2013}) or a depth-based version of kNN classification (\citealp{PaiVanB2012}), the dominant solution finds its source in the \emph{max-depth approach} of \cite{GhoCha2005B} that has later been refined by \cite{Lietal2012}. To the best of our knowledge, depth-based classification for directional data has not been considered in the literature. In this section, we show that the max-depth approach also applies for directional data and that, in conjunction with the proposed distance-based depths, it provides classifiers on the hypersphere that dominate ASD/ATD-based ones and that can be applied in higher dimensions as well. 

Consider the spherical classification problem where independent random samples~$W_{1i}$, $i=1,\ldots,n_1$ and~$W_{2i}$, $i=1,\ldots,n_2$, respectively, come from distributions~$H_1$ and~$H_2$ on~$\mathcal{S}^{q-1}$, and one is given the task to classify a point~$w(\in\mathcal{S}^{q-1})$ as arising from~$H_1$ (``population~1") or from~$H_2$ (``population~2"). Denoting as~$H_{\ell n_\ell}$ the empirical distribution associated with~$W_{\ell i}$, $i=1,\ldots,n_\ell$ ($\ell=1,2$),
the max-depth classifier associated with a depth~$D$ classifies~$w$ into population~1 if~$D(w,H_{1n_1})>D(w,H_{2n_2})$, and population~2 otherwise; if~$D(w,H_{1n_1})=D(w,H_{2n_2})$, then the classification decision is based on the flip of a fair coin.    

\begin{figure}[h!] 
\captionsetup{font=scriptsize}
\begin{center}  
\makebox[\textwidth]{%
\includegraphics[width=\textwidth]{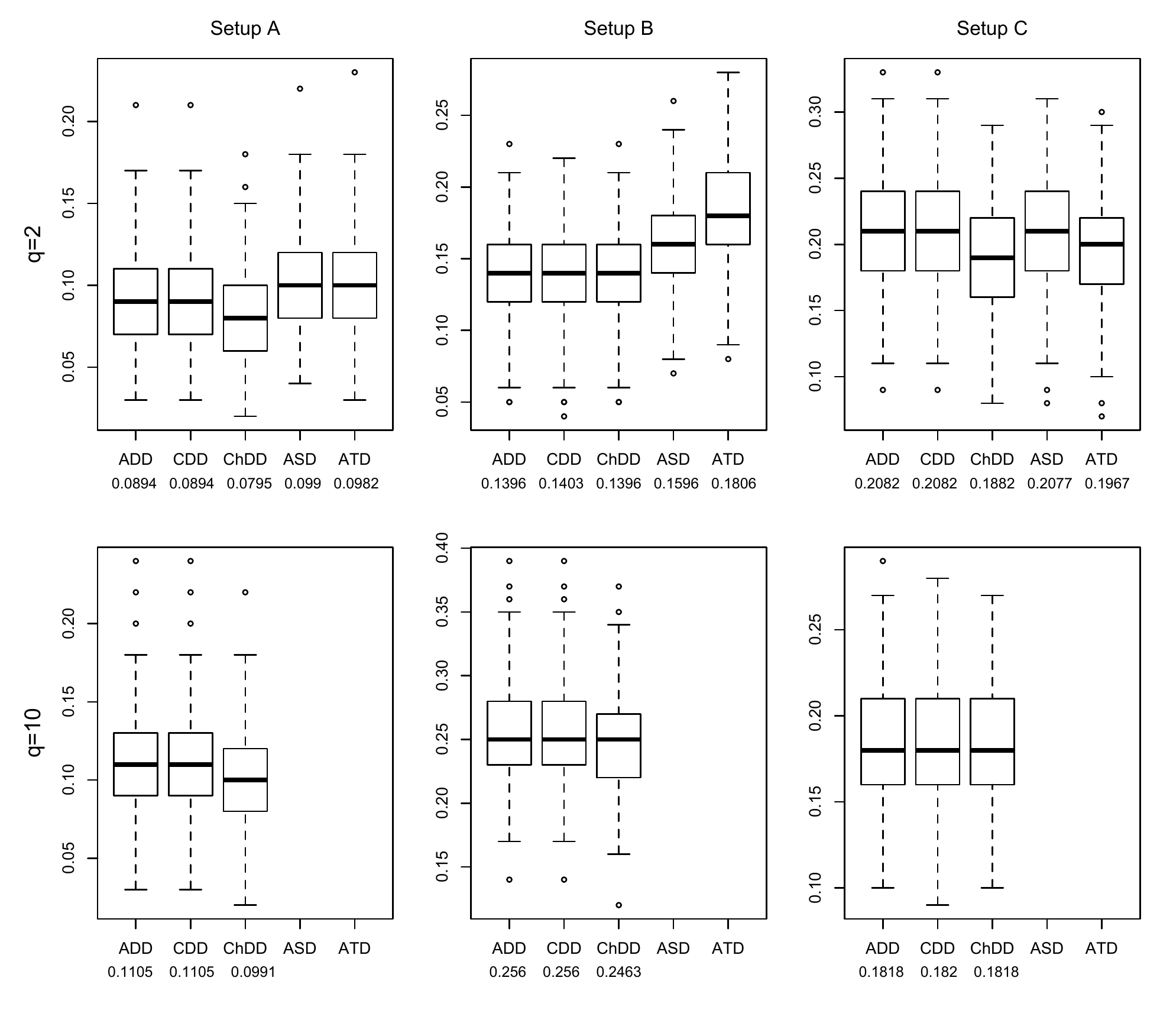}
}
\end{center}
\vspace{-5mm}
\caption{Boxplots, for~$q\in\{2,10\}$, of the misclassification rates~$p_m(D)$, $m=1,\ldots,M$, obtained from $M=250$ independent replications in three different distributional setups (see Section~\ref{secSimusub2} for details), for the max-depth classifiers associated with the ADD, CDD, ChDD, ASD and ATD (due to computational issues, the ASD and ATD were considered for dimension~$q=2$ only). In each case, the corresponding mean misclassification rate~$p(D) = (1/M) \sum_{m=1}^M p_{m}(D)$ is provided.} 
\label{Classifig}
\end{figure}

To investigate the finite-sample performances of such classifiers, we consider the Monte Carlo algorithm that was performed for dimensions~$q=2$ and~$q=10$. Denoting as~$e_j$ the $j$th vector in the canonical basis of~$\R^q$ and using the notations~$H^{\rm vMF}_{\alpha,\kappa}$ and~$H^{\rm vMF}_{q,\theta_1,\kappa}$ from Sections~\ref{secillu} and~\ref{sec:DistrProperties}, respectively, we considered the following three distributional setups:
\begin{itemize}
	\item Setup~A involves the vMF distributions~$H_1=H^{\rm vMF}_{\frac{\pi}{4},5}$ and~$H_2=H^{\rm vMF}_{\frac{3\pi}{4},5}$ for~$q=2$, and~$H_1=H^{\rm vMF}_{q,e_1,5}$ and~$H_2=H^{\rm vMF}_{q,e_q,5}$ for~$q=10$; Setup~A therefore involves distributions differing through the modal location only.
	\item In Setup~B,~$H_1=H^{\rm vMF}_{\frac{\pi}{3},2}$ and~$H_2=H^{\rm vMF}_{\frac{2\pi}{3},5}$ for~$q=2$, and~$H_1=H^{\rm vMF}_{q,e_q,2}$ and~$H_2=H^{\rm vMF}_{q,(\cos \frac{\pi}{6})e_{q-1}+(\sin \frac{\pi}{6})e_q,5}$ for~$q=10$; in this setup, distributions differ through location and concentration. 
	\item Setup~C involves discrimination between the vMF distribution~$H_1=H^{\rm vMF}_{\frac{3\pi}{4},4}$ and the mixture distribution~$H_2=\frac{1}{2}H^{\rm vMF}_{0,4}+\frac{1}{2}H^{\rm vMF}_{\frac{\pi}{2},4}$ for~$q=2$, and~$H_1=H^{\rm vMF}_{q,(\cos \frac{7\pi}{4})e_{q-1}+(\sin \frac{7\pi}{4})e_q,4}$ and~$H_2=\frac{1}{2}H^{\rm vMF}_{q,e_{q-1},4}+\frac{1}{2}H^{\rm vMF}_{q,e_q,4}$ for~$q=10$. 
\end{itemize}  
For each setup and each~$q$, we generated $M=250$ independent training samples of size~$n_{\rm train}=200$ and test samples of size~$n_{\rm test}=100$ by sampling randomly from~$\frac{1}{2}H_1+\frac{1}{2}H_2$. In replication~$m \in \{1,\ldots,250\}$, this associates with any depth~$D$ on~$\mathcal{S}^{q-1}$ the misclassification rate~$p_m(D)=N_m(D)/n_{\rm test}$, where~$N_m(D)$ is the number of observations in the $m$th test sample that were misclassified by the max-depth classifier associated with~$D$ when based on the $m$th training sample. Figure~\ref{Classifig} provides the boxplots, for several depths~$D$, of the resulting~$M=250$ misclassification rates. As in Section~\ref{secSimusub1}, the depths considered are the ADD, CDD, ChDD, ASD and ATD; again, computational issues prevented to consider the ASD and ATD in dimension~$q=10$.

Results indicate that distance-based depth classifiers dominate in most cases their counterparts based on the ASD/ATD. It is only in Setup~$C$ that the ASD/ATD classifiers seem to slightly improve over the ADD and CDD classifiers. In all cases, the classifier based on the ChDD is the best classifier. Most importantly, in higher dimensions, the computational burden for the ASD/ATD is such that only the distance-based depth classifiers can be used.


\section{Discussion}
\label{secfinal}

In the Euclidean multivariate setup, statistical depth has allowed to tackle in a nonparametric and robust way diverse problems, including location/scatter estimation, two-sample hypothesis testing, supervised classification, etc. While depths in the spherical setup, such as the ASD and ATD, were proposed more than two decades ago, the concept has not made its way to applications. Arguably, the reasons are that these depths are, even for moderate dimensions, very computationally intensive and that it is challenging to derive their asymptotic properties.  

The class of distance-based depths for directional data defined in this work clearly improve on this. These depths were showed to be computable in higher dimensions, and asymptotic results can be obtained by using standard $M$-estimation techniques. For small dimensions, where distance-based depths as well as the ASD/ATD can be evaluated, we showed through simulations that inference procedures based on the former compete equally or even dominate those based on the latter. In high dimensions, only distance-based depths can be used for directional data, which makes them of potential interest for applications involving high-dimensional spherical problems, such as those encountered in magnetic resonance, gene expression, or text mining; see, among others, \cite{Dry2005}, \cite{banerjee2003generative}, and \cite{Banetal2005}. 

Perspectives for future research are rich and diverse. Obviously, it would be of interest to investigate how distance-based depths can tackle the problems considered in the aforementioned high-dimensional applications. More generally, irrespective of the dimension, it would be desirable to develop depth-based inference procedures in various setups, including two-sample hypothesis testing and supervised classification. Finally, the present work also raised some theoretical questions of interest. For instance, in dimension~$q=2$, the arc distance depth is constant if and only if the underlying distribution~$H$ is antipodal, whereas the cosine distance depth is constant if and only if~$H$ has zero mean. In view of this, it is natural to wonder what property of~$H$ is characterized by constancy of the chord distance depth. The question can be raised on the circle with $q=2$ or for a general dimension~$q>2$. Such characterization results are of interest since they obviously provide the basis for universally consistent tests of the corresponding properties.


\appendix


As announced in Section~\ref{sec:StructProperties}, we prove the following result for the sake of completeness.

\begin{prop}
\label{prodistinv}
	Let~$d$ be a rotation-invariant distance on~$\mathcal{S}^{q-1}$. Then there exists a function~$\delta:[-1,1]\to\R^+$ such that~$d(\theta,\psi)=\delta(\theta'\psi)$.
\end{prop}

Proof of Proposition~\ref{prodistinv}.
For any~$\theta,\psi\in\mathcal{S}^{p-1}$, let~$\psi_\theta=(\psi-(\psi'\theta)\theta)/\|\psi-(\psi'\theta)\theta\|$ and denote as $\Gamma_{\theta,\psi}$ an arbitrary $q\times (q-2)$ matrix such that~$O_{\theta,\psi}=(\theta \vdots \psi_\theta \vdots \Gamma_{\theta,\psi})$ is orthogonal (if~$q=2$, then we simply consider~$O_{\theta,\psi}=(\theta \vdots \psi_\theta)$). Since~$d$ is rotation-invariant, we have~$d(\theta,\psi)=d(O_{\theta,\psi}'\theta,O_{\theta,\psi}'\psi)=d(e_1,O_{\theta,\psi}'\psi)$, where~$e_1$ stands for the first vector of the canonical basis of~$\R^q$. The result then follows from the fact that 
$
O_{\theta,\psi}'\psi
=
(
\theta'\psi , (1-(\theta'\psi)^2)^{1/2},0,\ldots,0
)'
$
depends on~$\theta$ and~$\psi$ through~$\theta'\psi$ only.
\cqfd
\vspace{3mm}


Proof of Theorem~\ref{thmrotainv}.
Using the notation introduced in the theorem, we have that $D_{d_{\delta}}(O\theta, H_O)
=\delta(-1)-E_{H_O}[\delta((O\theta)'W)]
=\delta(-1)-E_{H}[\delta((O\theta)'OW)]
=\delta(-1)
\linebreak
-E_{H}[\delta(\theta'W)]
=D_{d_{\delta}}(\theta, H)$.
\cqfd
\vspace{3mm}


Proof of Theorem~\ref{thcontinuity}.
(i) Since the function~$w\mapsto d(\theta,w)$ is continuous in~$w$ for any~$\theta\in\mathcal{S}^{q-1}$ and is bounded, uniformly in~$\theta$, by the integrable function~$w\mapsto d^{\sup}$, the continuity of
$$ 
\theta\mapsto D_{d}(\theta, H)
=
d^{\sup}
-
\int_{\mathcal{S}^{q-1}}
d(\theta,w)
\,
dH(w)
$$
results from Corollary 2.8.7(i) in \cite{Bog2007}. (ii) The result follows from the fact that  a continuous function on a compact domain attains its maximal value.
\cqfd
\vspace{3mm}


Proof of Theorem~\ref{thFishconsist}.
Since the distribution~$H$ is rotationally symmetric about~$\theta_0$, Theorem~\ref{thmrotainv} implies that~$D_{d_{\delta}}(\theta, H)$ depends on~$\theta$ only through~$\theta'\theta_0$. Consider then an arbitrary geodesic path~$t\mapsto \theta_t$ from~$\theta_0$ to~$\theta_1=-\theta_0$. The monotonicity assumption on~$h$ readily implies that, for any~$s\in[-1,1]$, the function~$t\mapsto P_{H}[ \theta_t'W \geq s]$ is monotone strictly decreasing. 
Since
\begin{eqnarray*}
E_H[\delta(\theta_t' W)]
&=&
\int_{0}^{\delta(-1)} z\, \frac{d}{dz} P_{H}[ \delta(\theta_t'W) \leq z] \,dz
\\[2mm]
&=&
\delta(-1) - \int_{0}^{\delta(-1)} P_{H}[ \delta(\theta_t'W) \leq z] \,dz
\\[2mm]
&=&
\delta(-1) - \int_{0}^{\delta(-1)} P_{H}[ \theta_t'W \geq \delta^{-1}(z)] \,dz
,
\end{eqnarray*}
it follows that
\begin{equation}
	\label{ahahah}
D_{d_{\delta}}(\theta_t, H)
=
\delta(-1)-E_H[\delta(\theta_t' W)]
=
\int_{0}^{\delta(-1)} P_{H}[ \theta_t'W \geq \delta^{-1}(z)] \,dz
\end{equation}
is strictly decreasing in~$t$. This establishes the result.
\cqfd
\vspace{3mm}


Proof of Theorem~\ref{propconcentr}.
First note that for any~$s$, 
\begin{equation}
	\label{tqhz1}
P_{H_\kappa}[\theta_0'W\geq s]
=
\frac{ 
\int_{s}^1
(1-v^2)^{(q-3)/2}
h(\kappa v) 
\,
dv
}
{\int_{-1}^1
(1-v^2)^{(q-3)/2}
h(\kappa v) 
\,
dv
}
\end{equation}
(see, e.g., \citealp{PaiVer17b}), which provides
\begin{equation}
	\label{tqhz2}
\frac{P_{H_\kappa}[\theta_0'W\geq s]}{1-P_{H_\kappa}[\theta_0'W\geq s]}
=
\frac{ 
\int_{s}^1
(1-v^2)^{(q-3)/2}
h(\kappa v) 
\,
dv
}
{\int_{-1}^s
(1-v^2)^{(q-3)/2}
h(\kappa v) 
\,
dv
}
\cdot
\end{equation}
Differentiation with respect to~$\kappa$ yields
\begin{eqnarray*}
\lefteqn{
\frac{d}{ds}
\frac{P_{H_\kappa}[\theta_0'W\geq s]}{1-P_{H_\kappa}[\theta_0'W\geq s]}
}
\\[2mm]
& & 
\hspace{3mm} 
=
\frac{ 
\int_{s}^1
\int_{-1}^s
[
v
\dot h(\kappa v) 
h(\kappa u) 
-
u
\dot h(\kappa u) 
h(\kappa v) 
]
((1-u^2)(1-v^2))^{(q-3)/2}
\,
du
dv
}
{
(
\int_{-1}^s
(1-v^2)^{(q-3)/2}
h(\kappa v) 
\,
dv
)^2
}
\cdot
\end{eqnarray*}
Since~$t \mapsto t\,\frac{d}{dt}\log h(t)=t \dot{h}(t)/h(t)$ is strictly increasing, this derivative is strictly positive at any~$\kappa$, so that the lefthand side of~(\ref{tqhz2}), hence also that of~(\ref{tqhz1}),  is a monotone strictly increasing function of~$\kappa$. The result then follows from the identity
$
D_{d_{\delta}}(\theta_0, H_\kappa)
=
\int_{0}^{\delta(-1)} P_{H_\kappa}[ \theta_0'W \geq \delta^{-1}(z)] \,dz
$;
see~(\ref{ahahah}). 
\cqfd
\vspace{3mm}


Proof of Theorem~\ref{thskewsym}.
(i) The anti-symmetry of~$\delta(\cdot)$ readily yields
$
D_{d_{\delta}}(-\theta, H) 
+
D_{d_{\delta}}(\theta, H) 
=
2 \delta(-1) - E_{H}[d_\delta(-\theta,W)+d_\delta(\theta,W)]
=
2 \delta(-1) - E_{H}[\delta(-\theta' W)+\delta(\theta' W)]
\linebreak
=
\delta(-1) 
$,	
which establishes the result. (ii) Ad absurdum, assume that~$-\theta_0$ does not have minimal depth, so that there exists~$\theta_1\in\mathcal{S}^{q-1}$ with~$D_{d_{\delta}}(\theta_1, H)<D_{d_{\delta}}(-\theta_0, H)$. Then Part~(i) of the result implies that~$D_{d_{\delta}}(-\theta_1, H)>D_{d_{\delta}}(\theta_0, H)$, which contradicts the fact that~$\theta_0$ has maximal depth.
\cqfd
\vspace{3mm}


Proof of Theorem~\ref{thunifconsistency}.
The result directly follows from Theorem~16(a) in \cite{Fer1996}.
\cqfd
\vspace{3mm}


Proof of Theorem~\ref{thasnormD}.
The result trivially follows from applying the central limit theorem to the expression~$
\sqrt{n}(D_{d}(\theta, H_n)-D_{d}(\theta, H))
=
-n^{-1/2} \sum_{i=1}^n (d(\theta,W_i)- {\rm E}_H[d(\theta,W)])
$.
\cqfd
\vspace{3mm}


Proof of Theorem~\ref{thdeepestconsistency}.
In view of Theorem~\ref{thunifconsistency}, the result is a corollary of Theorem~2.12 and Lemma~14.3 in \cite{Kos2008}.
\cqfd
\vspace{3mm}


Proof of Theorem~\ref{thmbdp}.
From Lemma~2.3 in \cite{str1985}, we obtain that, for any~$\theta\in\mathcal{S}^{q-1}$, 
$
| D_{d}(\theta,H_\varepsilon)-D_{d}(\theta,H) |
=
\varepsilon (E_G[d(\theta,W)]-E_H[d(\theta,W)])
\leq 
\varepsilon d^{\rm sup} d_1(H,G)
,
$
where $d_1(H,G)$ denotes the variational distance between~$H$ and~$G$. 
Lemmas~2.4 and~2.5(i) in \cite{str1985} then yield that, still for any~$\theta\in\mathcal{S}^{q-1}$,
$
| D_{d}(\theta,H_\varepsilon)-D_{d}(\theta,H) |
\leq
\varepsilon d^{\rm sup}
.
$
The result readily follows.
\cqfd
\vspace{3mm}




\bibliographystyle{chicago}
\bibliography{ManuscriptRevised.bib}           
\vspace{3mm}

\end{document}